\documentclass[preprint, number,5p,twocolumn,sort&compress,times,fleqn]{elsarticle}
\usepackage{amsmath}
\usepackage{amssymb}
\usepackage{mathtools}
\usepackage[colorlinks = true,
            linkcolor = blue,
            urlcolor  = blue,
            citecolor = blue,
            anchorcolor = blue]{hyperref}

\usepackage[utf8]{inputenc}
\usepackage{subfig}

\usepackage{csml}

\begin{document}
	
	\begin{frontmatter}
		
		\title{A three-grid high-order immersed finite element method for the analysis of CAD models}
		\author[1]{Eky Febrianto\corref{cor1}}
		\ead{eky.febrianto@glasgow.ac.uk}
		\cortext[cor1]{Corresponding author}
		
		\author[2]{Jakub {\v S}{\'\i}stek}
		\author[2]{Pavel K{\r u}s}
		\author[3]{Matija Kecman}
		\author[4]{Fehmi Cirak}  %Fehmi Cirak
		%\ead{f.cirak@eng.cam.ac.uk}
		
		\address[1]{Glasgow Computational Engineering Centre, University of Glasgow, University Avenue, Glasgow G12 8QQ, UK}
		\address[2]{Institute of Mathematics of the Czech Academy of Sciences, {\v Z}itn{\' a} 25, 115 67 Prague, Czech Republic}
		\address[3]{Epic Games, Inc. 620 Crossroads Blvd, Cary NC 27518, USA}
		\address[4]{Department of Engineering, University of Cambridge, Trumpington Street, Cambridge CB2 1PZ, UK}
		
		\begin{abstract}
The automated finite element analysis of complex CAD models using boundary-fitted meshes is rife with difficulties. Immersed finite element methods are intrinsically more robust but usually less accurate. In this work, we introduce an efficient, robust, high-order immersed finite element method for complex CAD models. Our approach relies on three adaptive structured grids: a geometry grid for representing the implicit geometry, a finite element grid for discretising physical fields and a quadrature grid for evaluating the finite element integrals. The geometry grid is a sparse VDB (Volumetric Dynamic B+ tree) grid that is highly refined close to physical domain boundaries. The finite element grid consists of a forest of octree grids distributed over several processors, and the quadrature grid in each finite element cell is an octree grid constructed in a bottom-up fashion. The resolution of the quadrature grid ensures that finite element integrals are evaluated with sufficient accuracy and that any sub-grid geometric features, like small holes or corners, are resolved up to a desired resolution. The conceptual simplicity and modularity of our approach make it possible to reuse open-source libraries, i.e. openVDB and p4est for implementing the geometry and finite element grids, respectively, and BDDCML for iteratively solving the discrete systems of equations in parallel using domain decomposition. We demonstrate the efficiency and robustness of the proposed approach by solving the Poisson equation on domains described by complex CAD models and discretised with tens of millions of degrees of freedom. The solution field is discretised using linear and quadratic Lagrange basis functions.
		\end{abstract}
		
		
	\end{frontmatter}
	
	%\tableofcontents
	
	%--------------------------------------------------------------------------------
\section{Introduction} 
\label{sec:introduction}
%--------------------------------------------------------------------------------
%
%--------------------------------------------------------------------------------
\subsection{Motivation and overview}
%--------------------------------------------------------------------------------
%
Immersed finite elements, also called embedded, unfitted, extended, cut-cell or shifted-boundary finite element methods, have unique advantages when applied to complex three-dimensional geometries~\cite{Hollig:2001aa,Hansbo:2002aa,freytag2006field,Parvizian2007,Ruberg:2010aa,Sanches:2011aa,Burman2015,Kamensky2015,main2018shifted,Noel2022,Schmidt2023,Li2024Complex}. By large, they can sidestep the challenges in automating the boundary-fitted finite element mesh generation from CAD models. The prevailing parametric CAD models based on BRep data structures and NURBS surfaces are generally rife with gaps, overlaps and self-intersections, which need to be repaired before mesh generation~\cite{shapiro2011geometric,marussig2017review,xiao2021delaunay}. In addition, CAD models contain small geometric details, like fillets, holes, etc., consideration of which would lead to unnecessarily large meshes so that they are typically manually removed, or defeatured~\cite{shapiro2011geometric}. Most of these challenges can be circumvented by present low-order immersed finite elements, but not so when high-order accuracy is desired. 

In immersed finite elements, the CAD model is discretised by embedding the domain in a structured background mesh, or grid, usually consisting of rectangular prisms, or cells. The boundary conditions are enforced approximately and the weak form of the governing equations is integrated only in the part of the cut-cells inside the domain. A polygonal approximation of the boundary surface will always lead to a low-order method irrespective of the polynomial order of the basis functions. For instance, according to standard error estimates for second-order elliptic boundary value problems to obtain an approximation of order ~$p+1$ using basis functions of degree~$p$, the boundary surface must be approximated with polynomials of degree~$p$~\cite{Ciarlet1972,Strang:2008aa}. That is, the cut-element faces adjacent to the boundary must be curved and be at least of degree~$p$. However, for complex geometries, generating geometrically and topologically valid curved elements is not trivial because of the mentioned problems with parametric CAD models. 

Alternatively, as pursued in this paper, the accuracy of the finite element solution can be increased by approximating the boundary surface as a polygonal mesh with a characteristic edge length much smaller than the cell size of the finite element grid. To this end, we introduce, in addition to the immersed \emph{finite element grid}, a \emph{geometry grid} for piecewise-linear implicit geometry representation and a \emph{quadrature grid} for evaluating the element integrals in the cut-cells. The three independent grids are all adaptive and have different resolutions. 

The conceptual simplicity of the proposed approach makes it easy to embed it within a parallel domain decomposition context and analyse large complex CAD models efficiently. Efficient solvers and parallel domain decomposition are essential for CAD geometries from engineering because their discretisation easily leads to systems of equations with several tens of thousands of elements and the associated computing times are below a few minutes. 
%
%\begin{figure*}[]
%	\centering
%	\subfloat[CAD model \label{fig:robotCAD}]{
%		\includegraphics[scale=0.33]{./figs/intro/robot_CAD-2}
%	}
%	\hspace{0.03\textwidth}
%	\subfloat[Polygonal STL mesh \label{fig:robotSTL}]{
%		\includegraphics[scale=0.13]{./figs/intro/roboSTL2}
%	}
%	\hspace{0.03\textwidth}
%	\subfloat[Signed distance function \label{fig:robotLevelset}]{
%		\includegraphics[scale=0.21]{./figs/intro/roboLvl}
%	}
%	\hspace{0.\textwidth}
%	\subfloat[Finite element solution \label{fig:robotDisp}]{
%		\includegraphics[scale=0.175]{./figs/intro/roboPotential2}
%	}
%	\caption[Immersed finite element analysis]{Robust higher-order accurate immersed finite element analysis of a CAD model. The BRep CAD model consisting of NURBS surface patches (a) is first exported from a CAD system as a sufficiently fine facetted STL mesh~(b). The respective implicit signed distance function~(c) is obtained by sampling the STL mesh.  The implicit signed distance function is stored on a grid restricted to a tight narrow band around the boundary. The finite element analysis is performed on a much coarser non-boundary-fitted grid using higher-order basis functions. The restriction of the finite element solution on the grid~(d) represents the sought solution.}
%	\label{fig:robot}
%\end{figure*}

\begin{figure*}[]
	\centering
		\subfloat[CAD model \label{fig:robotCAD}]{
			\includegraphics[scale=0.33]{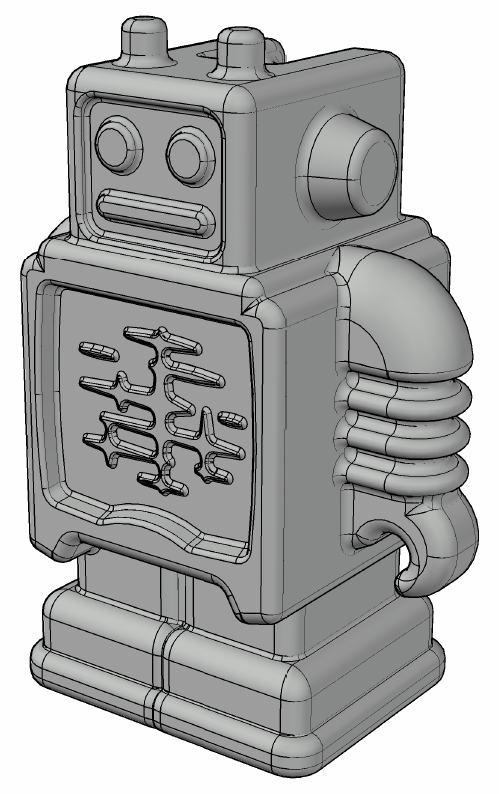}
		}
	\tikz[baseline=-\baselineskip]\draw[-stealth, line width=0.5mm] (1,2) --++  (1,0);
		\subfloat[Polygonal STL mesh \label{fig:robotSTL}]{
			\includegraphics[scale=0.13]{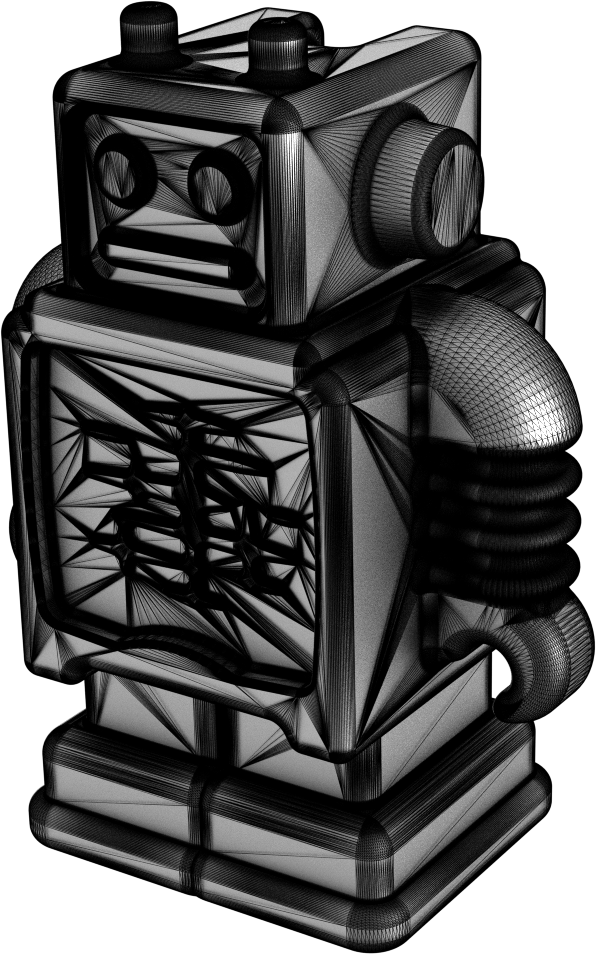}
		}
	\tikz[baseline=-\baselineskip]\draw[-stealth, line width=0.5mm] (1,2) -- ++ (1,0);
		\subfloat[Signed distance function \label{fig:robotLevelset}]{
			\includegraphics[scale=0.21]{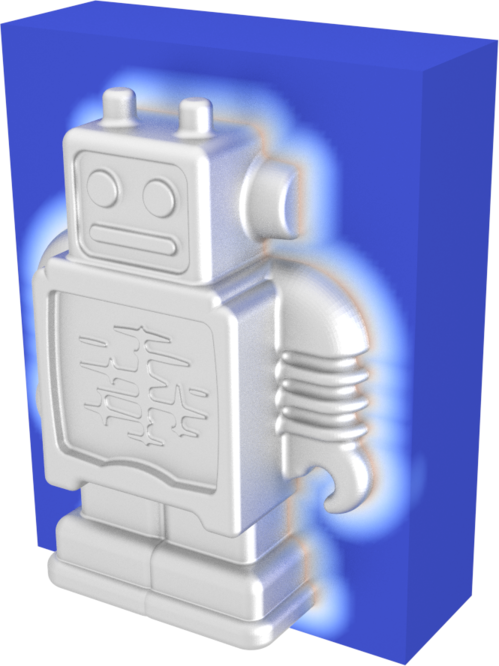}
		}
	\tikz[baseline=-\baselineskip]\draw[-stealth, line width=0.5mm] (1,2) -- ++ (1,0);
		\subfloat[FE solution \label{fig:robotDisp}]{
			\includegraphics[scale=0.175]{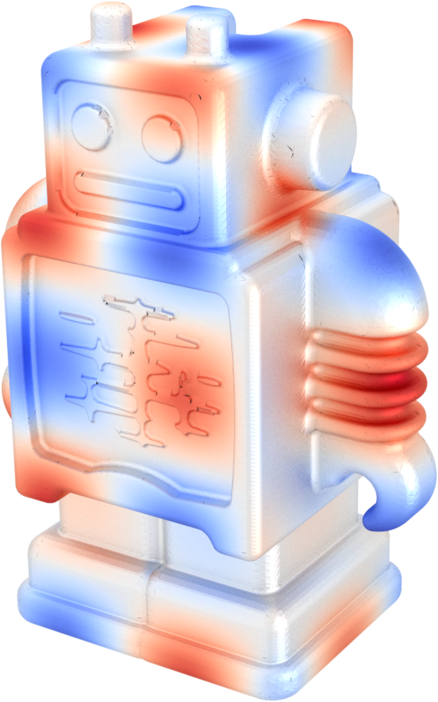}
		}
	\caption[Immersed finite element analysis]{Robust high-order accurate immersed finite element analysis of a CAD model. The BRep CAD model consisting of NURBS surface patches (a) is first exported from a CAD system as a sufficiently fine facetted STL mesh~(b). The respective implicit signed distance function~(c) is obtained by sampling the STL mesh.  The implicit signed distance function is stored on a grid restricted to a tight narrow band around the boundary. The finite element analysis is performed on a much coarser non-boundary-fitted grid using high-order basis functions. The restriction of the finite element solution on the grid~(d) represents the sought solution.}
	\label{fig:robot}
\end{figure*}

%
%--------------------------------------------------------------------------------
\subsection{Related research}
%--------------------------------------------------------------------------------
%
In an implicit geometry description, the domain is represented as a scalar-valued signed distance, or level set, function. The level set function is zero at the boundary, positive inside and negative outside of the domain. 
%Even when such functions do not encode the Euclidean distance to the boundary, they are usually called a signed distance function. 
Well-known advantages of implicit geometry representations include ease of robust Boolean operations and shape interrogation,  including the computation of ray-surface intersection, see e.g.~\cite{Ricci1973,farin2002handbook,museth2002level,patrikalakis2009shape,kambampati2021geometry}. 
%The approximate signed distance function for a given parametric CAD model, which will inevitably have problematic features, like gaps, overlaps and self-intersections. 
Although there are algebraic implicitisation techniques to determine the level set function of a single or a few NURBS patches~\cite{upreti2014algebraic,vaitheeswaran2021improved,xiao2019interrogation}, such methods are presently not scalable to complex CAD models. Alternatively,  the level set function can be obtained by sampling the CAD model and storing the determined distances on a structured geometry grid. However, sampling directly the CAD model requires distance computations or ray-surface intersections, leading to expensive and unstable nonlinear root-finding problems which would annihilate any advantages of an implicit geometry representation. 

In contrast, the sampling of polygonal surface meshes can be performed exceedingly robustly. Polygonal mesh models are prevalent in manufacturing and virtually all CAD systems can approximate a CAD model with an intersection-free polygonal mesh, that is an STL mesh, with a prescribed accuracy. An accurate representation of the level set function with a resolution comparable to the element size of the polygonal STL mesh is only required in the immediate vicinity of the boundary. Although octree data structures have been used to this end, they generally tend to have a large memory footprint and slow access times. 
%In a balanced octree, the refinement level between two neighbouring cells differs at most by one. 
%\ADDED{In a balanced octree, the branching factor is in each spatial dimension two, i.e., a parent cell yields eight sub-cells in 3D, which leads to a tall hierarchy.} 
An adaptive grid with approximately constant access time and a memory footprint that scales quadratically with the number of elements in the polygonal surface mesh is provided by the VDB (Volumetric Dynamic B+ tree)  data structure and its open source implementation OpenVDB~\cite{museth2013vdb,Museth2013}. 
In VDB the difference in the refinement level between two neighbouring cells is much higher than two, which is combined with other algorithmic features, making it ideal for processing and storing level set functions. 
%\ADDED{In VDB the branching factor is much higher than two yielding a shallow and wide tree with faster random access. This characteristic, combined with other algorithmic features make openVDB ideal for processing and storing level set functions.} 
OpenVDB also provides algorithms for efficient computation of the level set function by solving the Eikonal equation using the fast sweeping method~\cite{zhao2007parallel,zhao2005fast}.

The purpose of the finite element discretisation grid is to approximate the physical solution field and  usually has a different resolution requirement from the geometry grid for representing the level set function.   As per classical a-priori estimates, the ideal finite element cell size distribution depends on the smoothness properties of the solution field and the polynomial order of the basis functions used. A uniform cell size distribution is seldom ideal so that grid adaptivity is crucial. Furthermore, different from the geometry grid, abrupt changes in the finite element grid size distribution must be avoided because it is usually associated with artefacts in the approximate solution field. Hence, balanced octree grids are suitable as finite element discretisation grids. 

On the non-boundary-fitted octree finite element grid the boundary conditions are enforced by modifying either the variational weak form or the basis functions. According to historical work by Kantorovich and Krylov~\cite{kantorovichAndKrylov} in Ritz-like methods homogenous Dirichlet boundary conditions can be enforced by multiplying the global basis functions with a zero weight function at the boundary. Building on this idea, several approaches have been introduced for enforcing the boundary conditions by altering the basis functions~\cite{Rvachev1995,rvachev2000completeness,Hollig:2001aa,freytag2006field,Sanches:2011aa}. The required weight function is typically obtained from the signed distance function. Alternatively, the boundary conditions can be enforced by modifying the variational weak form, mainly using the Nitsche method~\cite{Nitsche:1971} and its variations~\cite{Embar2010,schillinger2016, Saye2017a, Saye2017b, Gulizzi2022}.  The finite element integrals in the cut-cells traversed by the boundary are evaluated only inside the domain.
Because finite element basis functions are locally supported, some basis functions close to cut-cells may have a negligible physically active support domain, leading to ill-conditioned system matrices. Several approaches have been proposed for cut-cell stabilisation, including basis function extension~\cite{Hollig:2001aa,Ruberg2012,Ruberg2016} and ghost-penalty method~\cite{burman2012fictitious}, see also the recent review~\cite{de2023stability}.

The evaluation of the finite element integrals in the cut-cells requires special care in high-order accurate immersed methods. Its accuracy depends on the chosen quadrature scheme, number of quadrature points, and the domain boundary approximation in the cut-cells. The number of quadrature points affects the overall efficiency of the immersed finite element method and must be kept low. The cut-cells are usually identified by evaluating the signed distance function at the cell vertices. This approach is sufficient when the geometry and finite element grids have the same resolution but can miss some of the cut-cells when the signed distance function is given in algebraic form or via a finer geometry grid. In both cases supersampling the signed distance function on the finite element grid can improve the correct classification of cells as cut-cells~\cite{luft2008geometrically}.  After identification, the cut-cells are decomposed into easy-to-integrate primitive shapes using marching tetrahedra partitioning~\cite{Ruberg2012,bandara2016shape}, octree partitioning~\cite{Parvizian2007} or a combination thereof~\cite{Kudela2016}.  For a discussion on the mentioned and other cut-cell integration approaches see~\cite{divi2020error}.  The collection of the primitive shapes, that is the integration cells, yields the integration grid. 

The quadrature approaches described so far are generally low order. In case of recursive octree partitioning high order accuracy can, in principle, be achieved by successively refining the integration cells traversed by the domain boundary. This leads however to a vast number of quadrature points and exceptionally inefficient approach.  An alternative approach to improve the boundary approximation is to introduce additional nodes on the edges and and faces and to push those to the boundary surface~\cite{Kudela2016,fries2017higher}. Unfortunately, such methods are very brittle because of the involved floating point operations, the presence of subgrid features, like small holes, and sharp corners and edges. They require a defeaturing step, like in boundary-fitted mesh generation, for engineering CAD models, annihilating one of the main advantages of immersed methods. It is worth bearing in mind, that the exceptional robustness of low order immersed methods can be attributed to the automatic defeaturing, or geometry filtering, by representing a CAD geometry on a coarse grid and the discarding of subgrid geometry details~\cite{freytag2006field,Sanches:2011aa}.

Parallel computation and grid adaptivity are essential for three-dimensional finite element analysis of large CAD models. The required computing memory and time for the analysis become quickly unwieldy, particularly when high-order basis functions are used. Adaptive refinement and coarsening of the finite element grid using a distributed octree data structure provides a mean to spread the grid and the computation over several processors and to choose a grid size distribution that best approximates the physical solution field~\cite{burstedde_p4est:_2011}. Compared to equivalent parallel boundary-fitted finite element implementations, an octree-based immersed finite element implementation leads to highly scalable domain partitioning and allows for easy dynamic load balancing, especially when adhering to balanced octree structures. There are efficient distributed octree implementations, most notably p4est~\cite{burstedde_p4est:_2011} for distributing the finite element grid over many processors. In such a distributed setting, the corresponding linear system of equations is typically solved using iterative Krylov subspace methods in combination with parallel preconditioners. The system matrix for the entire problem is never assembled. As parallel preconditioners, both algebraic multigrid~\cite{verdugo2019distributed} and domain decomposition~\cite{Badia-2018-RSD,jomo2019robust} have been applied in immersed methods. In~\cite{jomo2019robust} a standard single-level additive preconditioner and in~\cite{Badia-2018-RSD} a two-level balancing domain decomposition based on constraints~\cite{mandel2003convergence} are considered. The  possible ill-conditioning of the system matrices in immersed methods presents a challenge in applying preconditioners. However, this problem can be alleviated  by using the mentioned cut-cell stabilisation techniques.

%
%--------------------------------------------------------------------------------
\subsection{Proposed approach}
%--------------------------------------------------------------------------------
%	
The robustness, accuracy and efficiency of immersed methods effectively depend on the description of the domain geometry, the  treatment of the cut-cells, the resolution of the discretisation grid and the solution of the resulting linear systems of equations. In the case of large CAD models with geometric and topological faults and small geometric features,  an inevetiable trade-off is essential in satisfying the three competing objectives of robustness, accuracy and efficiency. To achieve this, we introduce three different adaptive grids for  representing the geometry, finite element discretisation and integration of element integrals. The resolution of each of the grids can be chosen in dependence of the required accuracy and available computing budget. As an additional benefit, the use of three different grids leads to a modular software architecture and makes it possible to reuse available open-source components, specifically openVDB~\cite{museth2013vdb}, p4est~\cite{burstedde_p4est:_2011} and BDDCML~\cite{Sousedik-2013-AMB}. 

The accuracy of the geometry representation is determined by the polygonal STL mesh with a user-prescribed precision exported from a CAD system. The respective signed distance function is computed at the nodes of the adaptive geometry grid with the minimum resolution length~$h_g$. The signed distance value within the cells is  obtained by linearly interpolating from the nodes. The length~$h_g$ is chosen in dependence of the smallest feature size appearing in the polygonal surface mesh. Any geometric details smaller than~$\approx h_g$ are automatically discarded by switching from the polygonal mesh to the implicit signed distance function. We use the openVDB library for computing and storing the implicit signed distance function in a narrow band of width~$\delta$ along the boundary surface. The narrow band size is chosen as a multiple of the minimum geometry resolution length~$h_g$. Inside the narrow band, the geometry cells have the length~$h_g$ and outside they are much coarser.  
%In the presented examples, we choose~$\delta = 3 h_g$. 

The adaptive finite element grid is constructed from a coarse base grid enveloping the polygonal surface mesh. The coarse base grid is distributed over the available processors and is refined until a desired cell size distribution is obtained while maintaining a refinement level difference of one or less between two neighbouring cells. The smallest resolution length~$h_f$ is usually larger than the geometry grid resolution~$h_g$ because the memory and computing requirements are dominated by the solution of the discretised equations. This is especially true when the domain is three-dimensional and high-order basis functions are used. In the included examples, we use Lagrange basis functions of polynomial orders~$p=1, \, 2$, and all are three-dimensional. Furthermore, all cut-cells are refined up to the finest resolution~$h_f$. We identify the cut-cells via supersampling the signed distance function with a resolution~$h_g$, within each finite element cell. This is necessary to detect small geometric features and sharp edges/corners that are easily missed by evaluating the signed distance function at the nodes of the finite element grid. We use the open source p4est (forest-of-trees) to refine and distribute the octree grid over all processors. Furthermore, we solve the respective distributed linear system of equations using the  adaptive-multilevel BDDC (balancing domain decomposition based on constraints library) BDDCML. 

The quadrature grid is constructed only in the cut-cells and has a resolution~$h_q \geq h_g$. It ensures that the finite element integrals are evaluated correctly in the presence of curved boundary surfaces and geometric features much smaller than the finite element cell size. We adopt a bottom-up octree strategy by starting from a coarse quadrature grid and successively merging cells that are not cut by the domain boundary. We implement fast cell traversal via the Morton code, or Z-curve. This bottom-up construction is crucial  for robust and accurate recovery of quadrature cells cut by the boundary. The leaf cells of the quadrature grid are then tessellated using the marching tetrahedra algorithm. 

%
%--------------------------------------------------------------------------------
\subsection{Overview of the paper}
%--------------------------------------------------------------------------------
%	
This paper is divided into five sections. In Section~\ref{sec:geometry}, we begin by providing a brief review of the implicit description of geometry and its discrete sampling over the geometry grid. We then introduce the FE grid and discuss its construction and classification in relation to geometrical and physical features. In Section~\ref{sec:immersed},  we present the immersed finite element scheme which utilises the three grids. We specifically address the treatment of cut-cells, including the construction of the bottom-up quadrature grid and tetrahedralisation. In Section~\ref{sec:parallel} we delve into the numerical strategies that enable the proposed immersed FE analysis to be computed on parallel processors. Finally, in Section~\ref{sec:examples}, we demonstrate the optimal convergence of the developed approach and its strong and weak scalability. We demonstrate the robustness of our approach through several examples, showing its ability to analyse various complex geometries. % See Figure~\ref{fig:all-geometry} for a visual representation.

	%--------------------------------------------------------------------------------
\section{Geometry and finite element grids}
\label{sec:geometry}
%--------------------------------------------------------------------------------

%-------------------------------------------------------------------------------
\subsection{Geometry grid} 
\label{sec:implicit}
%-------------------------------------------------------------------------------
%
We assume that the domain $\Omega \in \mathbb{R}^3$ with the boundary surface $\Gamma$ is given as a scalar valued implicit signed distance function $\phi(\vec x) : \Omega \rightarrow \mathbb R$ such that 
\begin{equation}\label{eq:sgndist1}
	\phi (\vec{x}) = 
	\begin{cases} 
		\phantom{-}  \dist( \vec{x}, \Gamma) \quad & \text{if } \vec x  \in \Omega \\ 
		\phantom{-}0 & \text {if } \vec{x} \in \Gamma \\
		-\dist(\vec{x},\Gamma) & \text{otherwise}\,,
	\end{cases}
\end{equation}
where  \mbox{$\dist(\vec x,\Gamma) = \min_{\vec y \in \Gamma} |\vec x - \vec y|$} denotes the shortest distance between the point $\vec x$ and the surface $\Gamma$. By definition, the signed distance function~$\phi(\vec x)$  is positive inside the domain, negative outside and the zeroth isosurface~$\phi^{-1}(0)$ corresponds to the boundary~$\Gamma$. 

\begin{figure*}[h!]
	\centering
	\subfloat[][Domain~$\Omega$ with the boundary~$\Gamma$ and the embedding domain~$\Omega_\Box$\label{fig:domain}]{
		\includegraphics[scale=0.16]{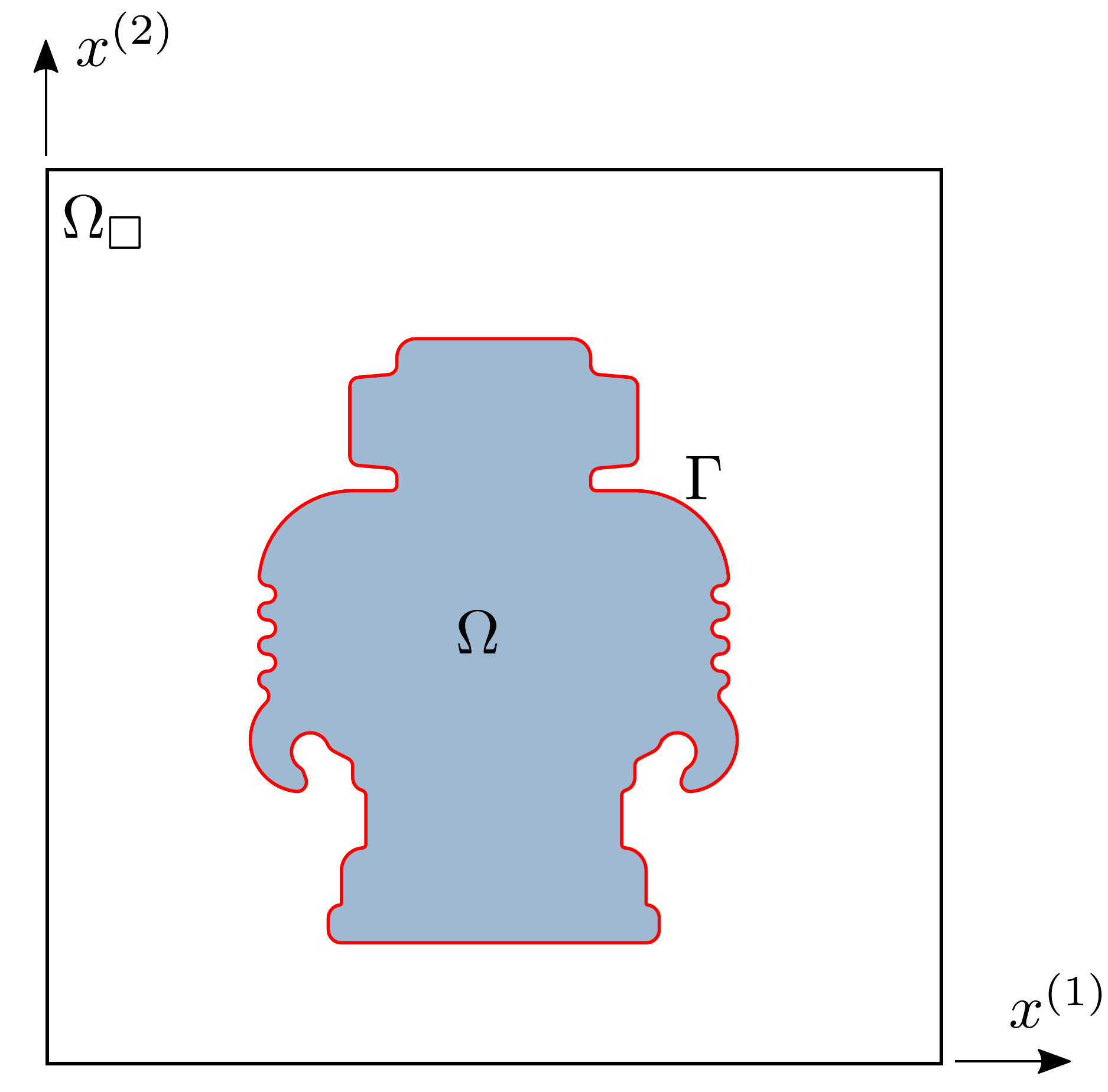}
	}
	\hfill
	\subfloat[][VDB grid with the smallest cell size $h_g$ of the geometry \label{fig:domainImplicitGrid}]{
		\includegraphics[scale=0.16]{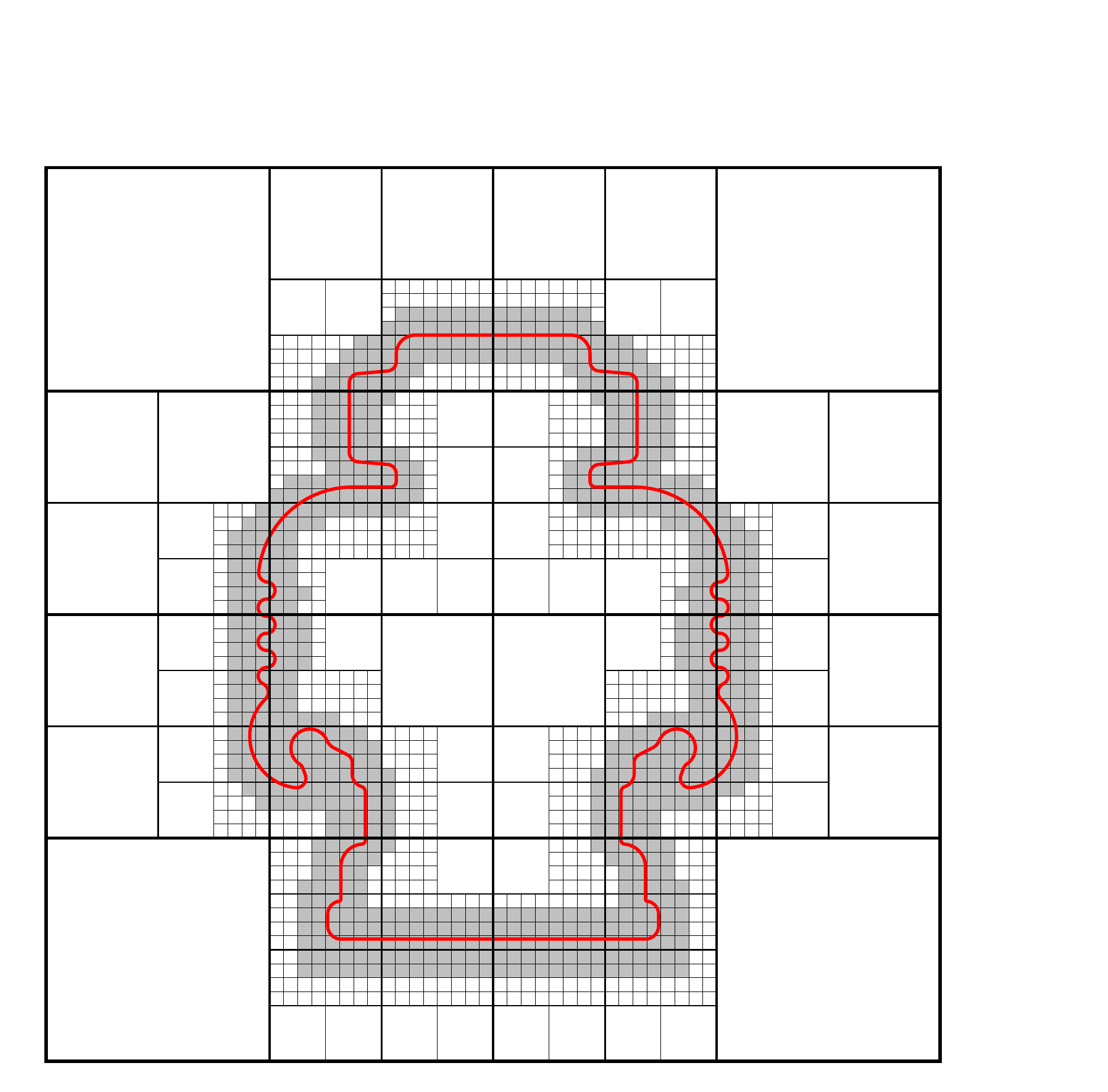}
	}
	\hfill
	\subfloat[][Approximate signed distance function~$\phi(\vec x)$ \label{fig:domainLevelSet}]{
	\includegraphics[scale=0.16]{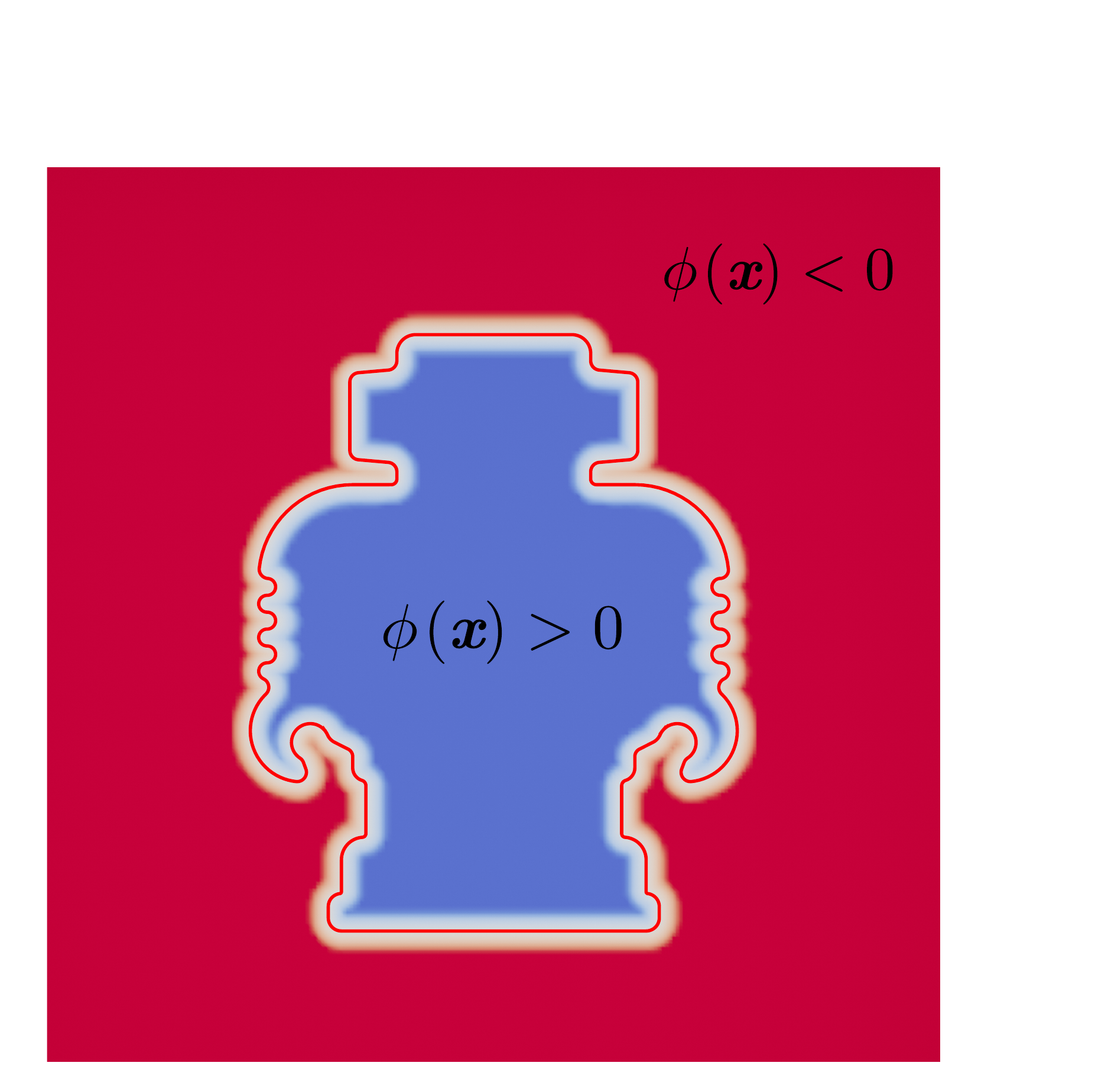}
	}
	\caption[]{Representative geometry setup.  \label{fig:domainDef}}
\end{figure*}

It is difficult to obtain the closed-form signed distance function~$\phi(\vec x)$ for complex CAD geometries especially when they consist of  freeform splines, see for example Figure~\ref{fig:robotCAD}. In that case, it is often expedient to consider a finely discretised polygonal surface mesh, i.e., STL mesh, approximating the exact CAD surface, see as an example Figure~\ref{fig:robotSTL} or its two-dimensional depiction in Figure~\ref{fig:domain}. Practically, all CAD systems are able to export STL mesh irrespective of the internal representation they use. 
%In that case, it is expedient to consider the discretised boundary~$\Gamma_t$ enclosing the domain~$\Omega_t$. We choose as~$\Gamma_t$ a sufficiently fine surface mesh, i.e. STL mesh, approximating the exact CAD surface~$\Gamma$, see as an example the two-dimensional illustration in~\ref{fig:domain}. Practically, all CAD systems are able to export STL mesh irrespective of the internal representation they use. 

To compute the approximate signed distance function~$\phi(\vec x)$, the respective domain~$\Omega$ is first embedded, or immersed, in a sufficiently large cuboid~\mbox{$\Omega_\Box \supset \Omega$}. We discretise the embedding domain~$\Omega_\Box$ with an adaptive grid consisting of cuboidal cells with minimum edge length~$h_g$. The signed distance function~$\phi(\vec x)$ is determined by first computing the distance of the grid nodes to the surface~$\Gamma$ and then linearly interpolating within the cells. The tree hierarchy of the geometry grid is established using the highly efficient VDB data structure as implemented in the open-source \mbox{openVDB} library~\cite{openvdbwebpage, museth2013vdb}.  The resulting geometry grid is highly refined within a narrow band of thickness~$2 \delta$ around the boundary, i.e., \mbox{$\{ \vec x \in \mathbb R^3 \mid \dist(\vec x,\Gamma) < \delta \}$}, see for example Figure~\ref{fig:domainImplicitGrid}, and has a very small memory footprint. In our computations, we choose~$\delta = 3~h_g$. The VDB data structure provides constant-time random access, insertion, and deletion, allowing fast scan conversion and evaluation of the signed distance function at any point in space. 
%In the two-dimensional robot illustration, the linearly interpolated signed distance function from the nodal samples is shown in Figure~\ref{fig:domainLevelSet}.

%
%--------------------------------------------------------------------------------
\subsection{Finite element grid}
\label{sec:femGrid}
%--------------------------------------------------------------------------------
%

\begin{figure*}[]
	\centering
	\subfloat[][Base finite element grid \label{fig:domainBaseGrid}]{
		\includegraphics[scale=0.18]{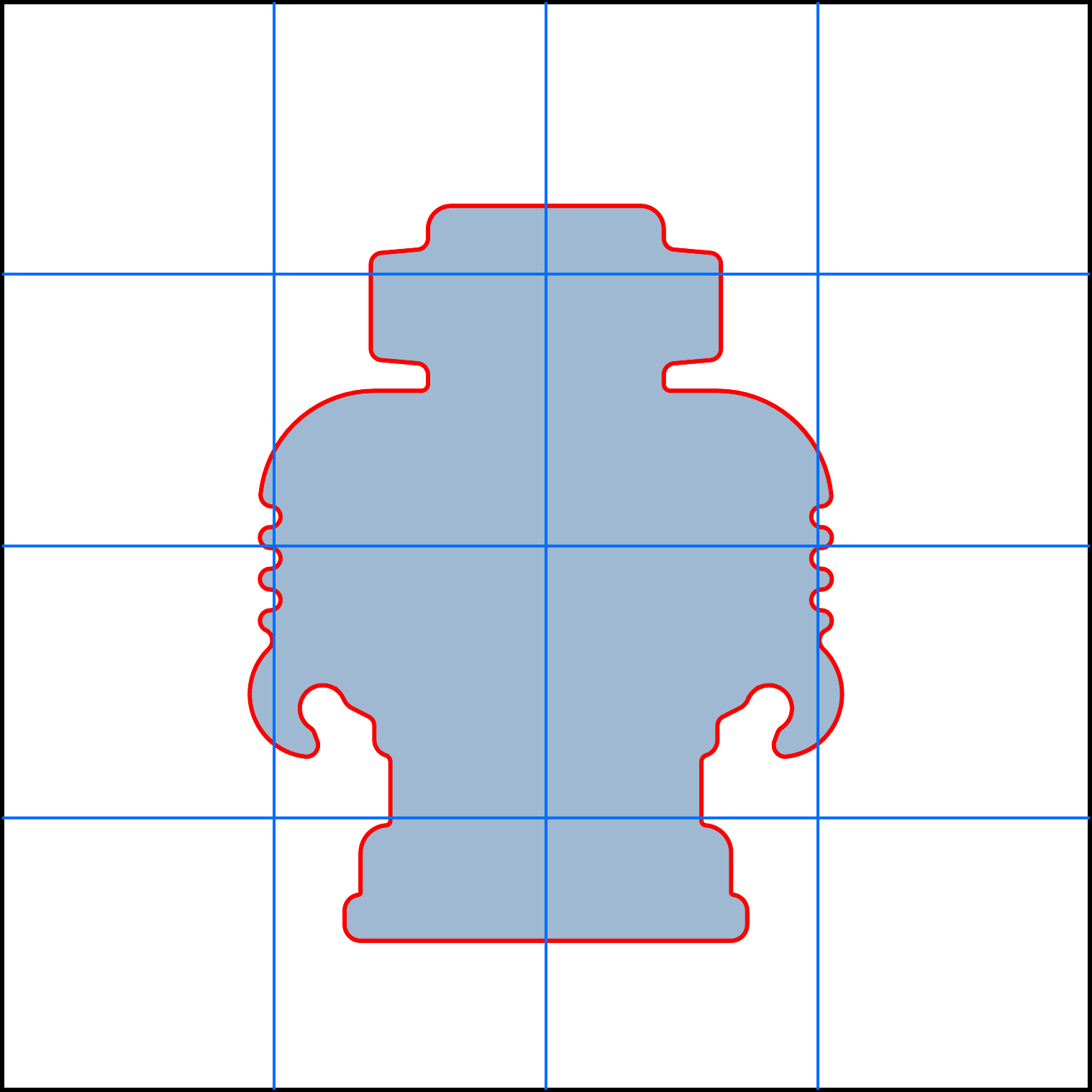}
	}
	\hspace{0.1\linewidth}
	\subfloat[][Finite element grid with smallest cell size $h_f$ \label{fig:domainFemGrid}]{
		\includegraphics[scale=0.18]{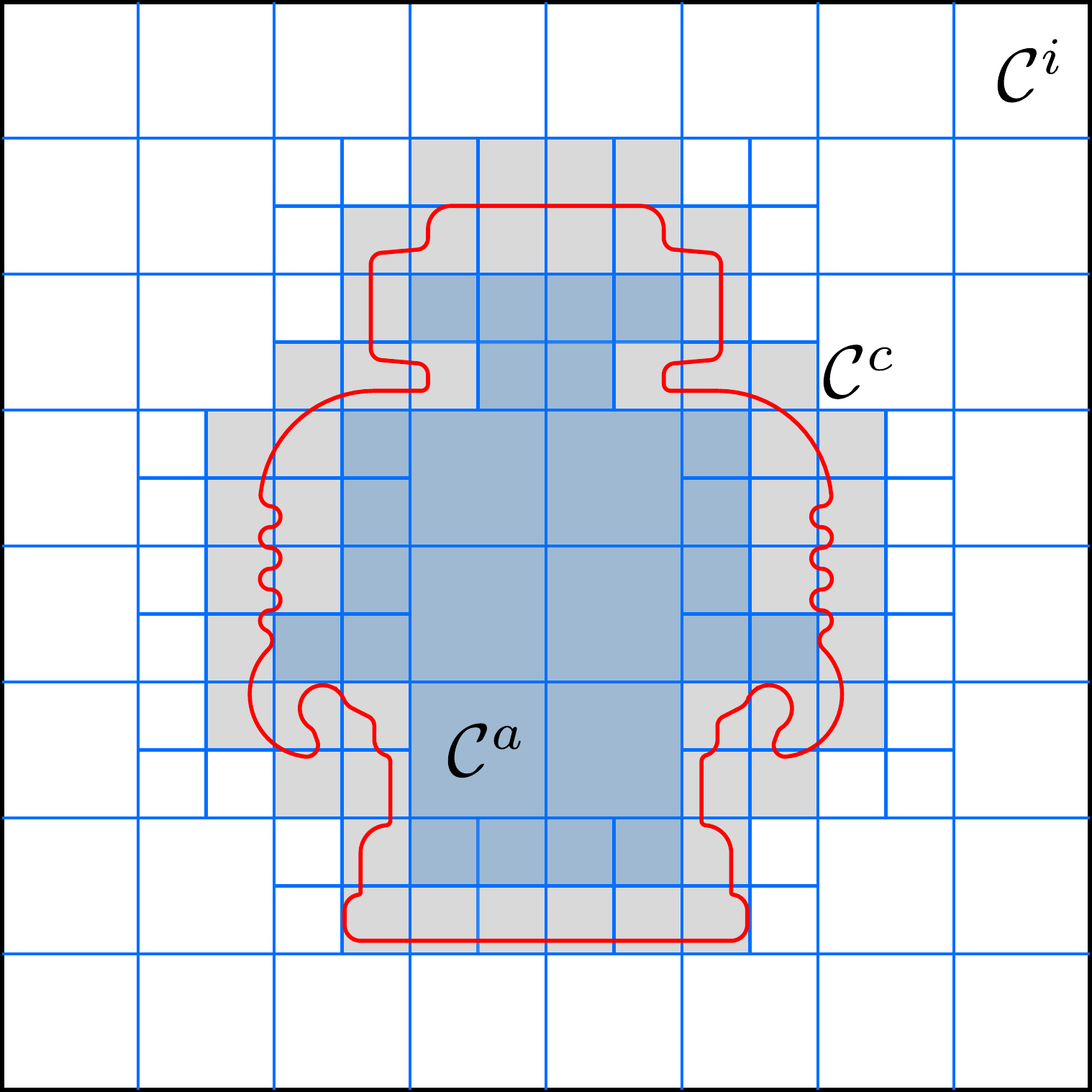}
	}
	\caption[]{A finite element grid obtained by octree refinement towards the domain boundary starting from a base grid. The domain boundary is given by the zeroth level set~$\phi^{-1}(0)$ of the signed distance function \label{fig:domainGrid}. The grid elements are categorised as active $\set C^a$, cut $\set C^c$ or inactive $\set C^i$.}
\end{figure*}

The finite element grid must ensure that the features of the physical solution field, for instance,  singularities and internal layers, are sufficiently captured. 
%In the immersed FE context, where the discretisation often do not coincide with the geometric boundary, the FE grid is also required to correctly identify the boundary presence. 
The construction of the FE grid starts with a very coarse base grid, see Figure~\ref{fig:domainBaseGrid}. The base grid is usualy obtained by uniformly subdividing the bounding box $\Omega_{\Box}$. To obtain the FE grid, we increase the resolution of the base grid by repeated octree refinement of selected cells, i.e. by subdividing cells into eight cells. The maximum and minimum cell size are often prescribed as constraints of the refinement. 

In boundary value problems it is usually necessary to refine the grid towards the domain boundaries because of the practical importance of the solution close to the boundaries.  As an example, we illustrate in Figure~\ref{fig:domainFemGrid} the refinement of the base grid towards the boundary utilising the FE grid cell predicates that will be introduced in Section~\ref{sec:cutCells}. When a cell is identified as a cut-cell, it is subdivided into eight sub-cells. The octree data structure representing the finite element grid is managed using the open-source library p4est~\cite{burstedde_p4est:_2011}.  As will be discussed in Section~\ref{sec:partitioning},  p4est also manages the partitioning of the octree into subprocesses and distributes them across several processors.  Different from the VDB data structure for the geometry, we constrain the refinement level between adjacent cells to have a 2:1 ratio, i.e., the refinement level between neighbouring cells may differ at most by one. 
%This restriction is sensible in case of physical fields without highly localised singularities, and \code{p4est} is able to enforce this requirement.

The resulting finite element grid comprises of the non-overlapping leaf cells~$ \set C = \{ \omega_i \}$ that decompose the bounding box, i.e., 
\begin{equation}\label{eqn:feCellDiscretisation}
	\Omega_\Box = \bigcup_i \, \omega_i \, .
\end{equation}
Each cell~$\omega_i$ represents a finite element, and the physical field is discretised using the basis functions associated with the cells. 

%-------------------------------------------------------------------------------
\subsection{Classification of finite element cells \label{sec:cutCells}}
%-------------------------------------------------------------------------------
%
In this section we describe the classification of FE grid cells pertinent to the adaptive refinement of the cells. It is assumed that the discretised signed distance function~$\phi(\vec x)$ can be evaluated at any point $\vec x \in \Omega_{\Box}$. 
%Clearly, the signed distance function~$\phi_g(\vec x)$  provides only an approximation to the domain~$\Omega_t$ given by the triangle surface mesh. At a particular refinement step, 
The finite element cells~$\set C$ are split into three disjoint sets
\begin{subequations} \label{eq:cellSets}
	\begin{align}
		\set C^a &= \{ \omega_{\vec i} \in \set C  | \min_{\vec x  \in \omega_{\vec i}} \phi (\vec x)   > 0 \} \, , \\
		\set C^i &= \{ \omega_{\vec i} \in \set C  | \max_{\vec x  \in \omega_{\vec i}} \phi (\vec x )   < 0 \} \, ,\\
		\set C^c & = \set C \setminus (  \set C^a \cup \set C^i )   \, . 
	\end{align}
\end{subequations}
The cells in the sets~$\set C^a$, $\set C^i$ and $\set C^c$ are referred to as the active, inactive and cut cells, see Figure \ref{fig:domainFemGrid}. 

\begin{figure*}[]
	\centering
	\subfloat[][\label{fig:nonOrdinaryFeaturesA}]{
		\includegraphics[scale=1]{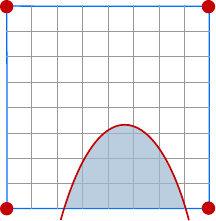}
	}
	\hspace{0.09\textwidth}
	\subfloat[][\label{fig:nonOrdinaryFeaturesB}]{
		\includegraphics[scale=1]{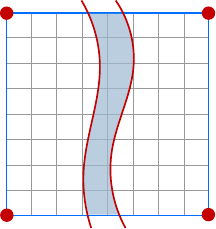}
	}
	\hspace{0.09\textwidth}
	\subfloat[][\label{fig:nonOrdinaryFeaturesC}]{
		\includegraphics[scale=1]{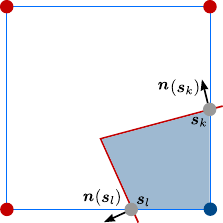}
	}
	\caption[]{Identification of the extraordinary cut finite element cells~$\omega_i \in \set C^{c_e} \subset \set C$. Blue dots correspond to $\phi(\vec x) > 0 $ and red dots to $\phi (\vec x) < 0 $. Supersampling of selected cells to detect intersections with small geometric details (left and centre) and cells with sharp geometric features (right). In (a) and (b) $h_f = 8 h_g$. }
	\label{fig:nonOrdinaryFeatures}
\end{figure*}

Although it is tempting to classify the finite element cells by evaluating the signed distance function~$\phi(\vec x)$ only at their corners, it usually leads to a crude finite element approximation of the domain~$\Omega$. As depicted in Figures~\ref{fig:nonOrdinaryFeaturesA} and~\ref{fig:nonOrdinaryFeaturesB},  such an approach may miss physically important small geometric details.  A more accurate finite element approximation can be obtained by supersampling the signed distance function~$\phi(\vec x)$ within the cells. The chosen sampling distance represents a low-pass filter for the geometry and  implies a resolution length for geometric details. The sampling distance can be equal or larger than the edge length~$h_g$ of the geometry grid. 
%In addition to this geometric filtering, another advantage of classifying the cells using~$\phi_g(\vec x)$ rather than directly~$\Omega_t$ is that it leads to a robust approach while not requiring intricate intersection computations. 

\begin{figure*}[]
	\centering
	\subfloat[][\label{fig:ordinaryFeaturesA}]{
		\includegraphics[scale=1]{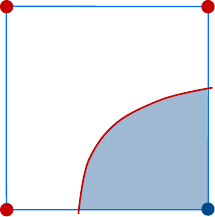}
	}
	\hspace{0.09\textwidth}
	\subfloat[][\label{fig:ordinaryFeaturesB}]{
		\includegraphics[scale=1]{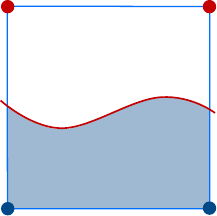}
	}
	\hspace{0.09\textwidth}
	\subfloat[][\label{fig:ordinaryFeaturesC}]{
		\includegraphics[scale=1]{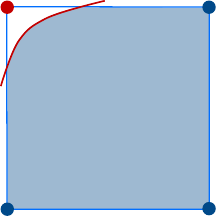}
	}
	\caption[]{Identification of the ordinary cut finite element cells~$\omega_i \in \set C^{c_o} \subset \set C$. Blue dots correspond to $\phi (\vec x) > 0 $ and red dots to $\phi (\vec x) < 0 $.}
	\label{fig:ordinaryFeatures}
\end{figure*}

Industrial CAD geometries usually contain sharp features in the form of creases and corners that are often physically essential, see Figure~\ref{fig:nonOrdinaryFeaturesC}. They can be identified by considering the change of the normal to the signed distance function~\cite{Kobbelt2001}. The respective unit normal~$\vec n(\vec x)$ is given by 
\begin{equation}
	\vec n (\vec x) = \frac{\nabla \phi(\vec x) }{ \| \nabla \phi (\vec x) \|} \, ,
\end{equation}
where~$\nabla$ denotes the gradient operator.

Consider a set of points in a cell $\omega_i$ denoted as~$\set S_i$, where we can evaluate $\phi$ and $\nabla \phi$. In our implementation the set~$\set S_i$ corresponds to grid points of a local uniform grid. If the spacing of this grid is the cell size itself, $\set S_i$ contains just the vertices of $\omega_i$. We can define an indicator of sharp features by
\begin{equation} \label{eq:sharpCrit}
	 \min_{\vec s_k  \in \set S_i, \, \vec s_l \in \set S_i} \vec n (\vec s_k ) \cdot \vec n (\vec s_l) < \cos(\theta) \, , 
\end{equation}
where~$\cos(\theta)$ is a user prescribed parameter chosen as~$\cos(\theta) = 0.3$ in the presented computations. 

In light of the mentioned observations, we employ the following algorithm to classify the finite element cells. 
\begin{itemize}
	\item[S1.] Define ${\set S_i}$ as the vertices of the finite element cell $\omega_i$, and
                   sample and store the signed distance values~$\phi(\vec s_j)$ and the normal vectors $\vec n(\vec s_j)$
                   for $\vec s_j \in {\set S_i}$.
	\item[S2.] Identify all the cut-cells $\set C^{c}$ such that 
	\begin{equation} \label{eq:ordinaryCutCrit}
		\set C^{c} = \{ \omega_{\vec i}  \mid  \,  \min_{\vec s_j  \in \set S_i} \phi (\vec s_j) \cdot \max_{\vec s_k  \in \set S_i} \phi (\vec s_k)  < 0 \}\, ,
	\end{equation}
	identify the subset of cut-cells $ \set C^{c_e} \subset \set C^{c} $ that contain a sharp feature using criterion~\eqref{eq:sharpCrit}, and define the subset of \emph{ordinary} cut cells $\set C^{c_o} := \set C^{c}  \backslash  \set C^{c_e}$, as depicted in Figure~\ref{fig:ordinaryFeatures}.
    Categorise the remaining cells into active $\set C^a$ and inactive $\set C^i$ considering the signed distance values~$\phi(\vec s_{j})$.
	\item[S3.] Collect the set of cells $\set P = \bigcup_i \set P_i$, where $\set P_i$ indicates the face neighbour of a cut-cell $\omega_i \in \set C^{c}$.  
	\item[S4.] For cells in $\set P \setminus \set C^{c}$, redefine $\set S_i$ by supersampling, see for example Figure~\ref{fig:nonOrdinaryFeatures}, 
        and repeat S1 and S2. 
        Cells identified as cut after supersampling are added to the set of extraordinary cells $\set C^{c_e}$.
        
\end{itemize}
The final set of the cut cells is given as $\set C^c = \set C^{c_o} \cup \set C^{c_e}$.

\nobreak 
	
	%
%-------------------------------------------------------------------------------
\section{Immersed finite element method}
\label{sec:immersed}
%-------------------------------------------------------------------------------

%-------------------------------------------------------------------------------
\subsection{Governing equations}
\label{sec:governingEqn}
%-------------------------------------------------------------------------------
%
As a representative partial differential equation, we consider on a domain~\mbox{$\Omega \in \mathbb R^3$}  with the boundary~\mbox{$\Gamma = \Gamma^D \cup \Gamma^N$} the Poisson equation 
\begin{equation}\label{eq:femPoissonStrongForm}
	\begin{aligned}
		-\nabla \cdot \nabla u &= f && \text{in $\Omega$} \, ,
		\\
		u  &  = \overline{u}  && \text{on $\Gamma^D$}  \, ,
		\\
		\vec{n}  \cdot \nabla{u}  &  = \overline{g}  && \text{on $\Gamma^N$} \, ,
	\end{aligned} \, 
\end{equation}
where~\mbox{$u, f: \Omega \rightarrow \mathbb R$} are the unknown solution and the prescribed forcing,  \mbox{$\overline{u}: \Gamma^D \rightarrow \mathbb R$} is the prescribed Dirichlet data,~$\overline{g}: \Gamma^N \rightarrow \mathbb R$ is the prescribed Neumann data, and~$\vec n:  \Gamma \rightarrow \mathbb R^3$  is the unit boundary normal.  The weak formulation of the Poisson equation can be stated according to Nitsche~\cite{Nitsche:1971} as follows: find $u \in H^1(\Omega)$ such that
\begin{equation}\label{eq:femWeakForm}
	\begin{aligned}
	\int_\Omega \nabla{u} \cdot \nabla{v}  \D \Omega  & +  \gamma \int_{\Gamma_D} (u  -\overline{u})  v  \dif \Gamma \\
	&= 
	\int_\Omega f  v \D \Omega + \int_{\Gamma_N} \overline{g}  v \dif \Gamma \\
	&+ 
	\int_{\Gamma_D} \left (   ( u -\overline{u} ) \vec{n}  \cdot \nabla v  + (\vec{n} \cdot \nabla u)\, v \right )  \dif \Gamma 
	\end{aligned}
\end{equation}
for all $v \in  H^1(\Omega)$. The space $ H^1(\Omega)$ is the standard Sobolev space  such that the test functions $v$ do not have to be zero on~$\Gamma_D$ and the Dirichlet boundary conditions are satisfied only weakly by the solution $u$. The stability parameter  $\gamma > 0$ can be chosen, for instance, according to~\cite{Embar2010,Ruberg2016} as $\gamma_0 / h_f $, where $h_f$ is the local finite element size.
This choice of $\gamma$ requires that the adverse effects of small cut cells is eliminated by other means, e.g., by basis function extrapolation as will be introduced in Section~\ref{sec:stabilisation}.

We discretise the solution field~$u$ and the test function~$v$ with Lagrange basis functions~$N_i(\vec x)$ that are defined on the non-boundary fitted finite element grid defined on the domain~$\Omega_\Box \supset \Omega$. The approximation of the solution~$u$ and the test function~$v$ over the finite element grid is given by
\begin{equation}
	\label{eq:interpolation}
		u^h  = \sum_{\vec{i}}  N_{\vec{i}} (\vec x  ) u_{\vec{i}} \, ,   \quad  v^h  = \sum_{\vec{i}}  N_{\vec{i}} (\vec x ) v_{\vec{i}} \,  . 
\end{equation}
After introducing both into the weak form~\eqref{eq:femWeakForm} all integrals are evaluated by iterating over the cells in the finite element grid. Only the active cells in~$\set C^a$ and the cut cells in~$ \set C^c$ have a non-zero contribution and are considered in finite element analysis. The integrals over active cells can be evaluated using standard tensor-product Gauss quadrature. The evaluation of the integrals over cut cells requires some care. As mentioned, according to standard a-priori error estimates, optimal convergence of the finite element solution is only guaranteed when the boundary geometry is approximated with polynomials of the same degree as the basis functions~$N_i(\vec x)$. 
%, see e.g.~\cite{Strang:2008aa}. 

Unfortunately, trying to reconstruct the boundary geometry from the zeroth isocontour of the discretised signed distance function, i.e.~$\phi^{-1}(0)$,  with higher than linear polynomials is as challenging as creating a conforming 3D finite element mesh. Therefore we use in the proposed approach only a piecewise linear reconstruction of the boundary on the cut cells in~$\set C^c$, as will be detailed in Section~\ref{sec:cutTesselate}. To increase the integration accuracy, we introduce a fine quadrature grid obtained by octree refinement of the finite element cut cells, which will be detailed in Section~\ref{sec:bottomUpCut}. The discretisation of the weak formulation~\eqref{eq:femWeakForm} yields after integration using the quadrature grid the discrete system of equations
\begin{equation} \label{eq:discEquilibrium}
	\vec A \vec u = \vec f \, ,
\end{equation}
where~$\vec A$ is the symmetric positive definite system matrix, $\vec u$ is the solution vector and $\vec f$ is the forcing vector.

%
%-------------------------------------------------------------------------------
\subsection{Quadrature grid and bottom-up cut cell refinement} 
\label{sec:bottomUpCut}
\begin{figure}[h!]
	\centering
	\subfloat[][Cut finite element cell \label{fig:sharpSubgridA}]{
		\includegraphics[scale=1]{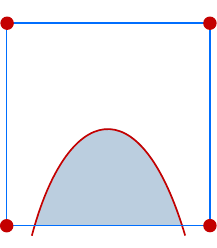}
	}
	\hspace{0.02\textwidth}	
	\subfloat[][Z-curve of the leaf cells ($r_Q=3$)\label{fig:sharpSubgridB}]{
		\includegraphics[scale=1]{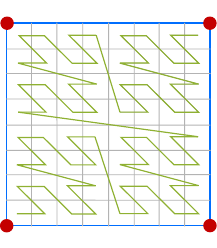}
	}\\
	\subfloat[][Quadrature cells \label{fig:sharpSubgridC}]{
		\includegraphics[scale=1]{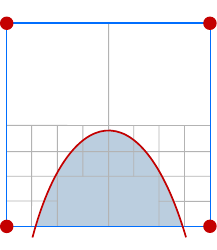}
	}
	\hspace{0.02\textwidth}	
	\subfloat[][Tetrahedralisation of cut quadrature cells\label{fig:sharpSubgridD}]{
		\includegraphics[scale=1]{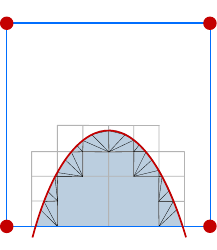}
	}
	\caption[]{Quadrature grid and bottom-up cut cell refinement of a cut finite element cell in~$\set C^c$ for quadrature and the tetrahedralisation of the cut quadrature cells. The red dots correspond to $\phi (\vec x) < 0 $.}   
	\label{fig:sharpSubgrid}
\end{figure}

%-------------------------------------------------------------------------------
%
The quadrature grid in the cut finite element cells in~$ \set C^c$ is obtained by~$r_Q$ steps of octree refinement.  This implies that $r_Q = h_f / h_q$, where $h_q$ is the smallest quadrature cell size. In our computations, we usually choose~$0 \le r_Q \le 4$. The refinement process yields for each cut cell~$8^{r_Q}$ leaf quadrature cells over which the finite element integrals have to be evaluated. As exemplified in Figure~\ref{fig:sharpSubgridA}, not all of the leaf quadrature cells are intersected by the domain boundary represented by $\phi^{-1}(0)$. Considering that the evaluation of the finite element integrals is computationally expensive (more so for high-order basis functions), it is desirable to have as few integration cells as possible. Hence, it is expedient to merge all the leaf cells not intersected by~$\phi^{-1}(0)$ into larger cuboidal cells. This does not harm the finite element convergence rate because the resulting larger cells are integrated with the same tensor product quadrature rule like the active elements in~$\set C^a$.  

We use a bottom-up octree construction to merge the leaf cells of the quadrature octree grid into larger integration cells. To this end, first, all the leaf quadrature cells are classified into active, inactive and cut by evaluating the signed distance value~$\phi_g(\vec x)$ at their corners. Subsequently, the bottom-up octree  is constructed using the Morton code, or Z-curve~\cite{Morton1966}, of the leaf quadrature cells, see Figure~\ref{fig:sharpSubgridB}. As the Z-curve is locality preserving, we can efficiently traverse it and apply a simple coarsening to build up the octree from the bottom up. The so obtained few larger cells can be integrated much more efficiently, see Figure~\ref{fig:sharpSubgridC}. For future reference, we denote the subcells of the finite element cut cell~$\omega_i \in \set C^c$ with~$\omega_{i, k}^Q$ such that 
\begin{equation} \label{eq:quadCells}
	\omega_i =  \bigcup_k   \omega_{i,k}^Q \, .
\end{equation}
In the following, the cells~$\omega_{i, k}^Q$ are referred to as quadrature cells.
Active quadrature cells are integrated using the tensor-product Gauss quadrature in a straightforward way. On the other hand, the cut quadrature cells first undergo the tetrahedralisation algorithm followed by integration introduced in Sections~\ref{sec:cutTesselate} and~\ref{sec:feIntegration}.

\begin{figure*}[]
	\centering
	\subfloat[][Octree refinement level~$r_Q=0$ \label{fig:sharpFeatureCracksA}]{
		\includegraphics[scale=0.55]{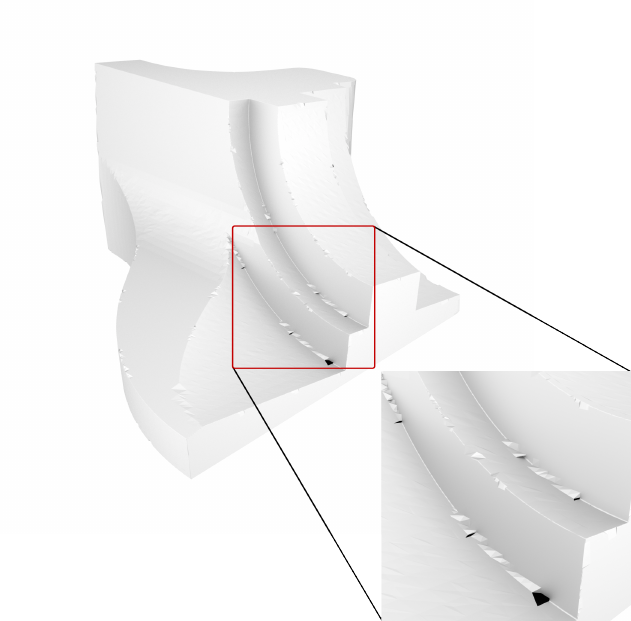}
	}
	\hspace{0.05\textwidth}	
	\subfloat[][Octree refinement level~$r_Q=3$ \label{fig:sharpFeatureCracksB}]{
		\includegraphics[scale=0.55]{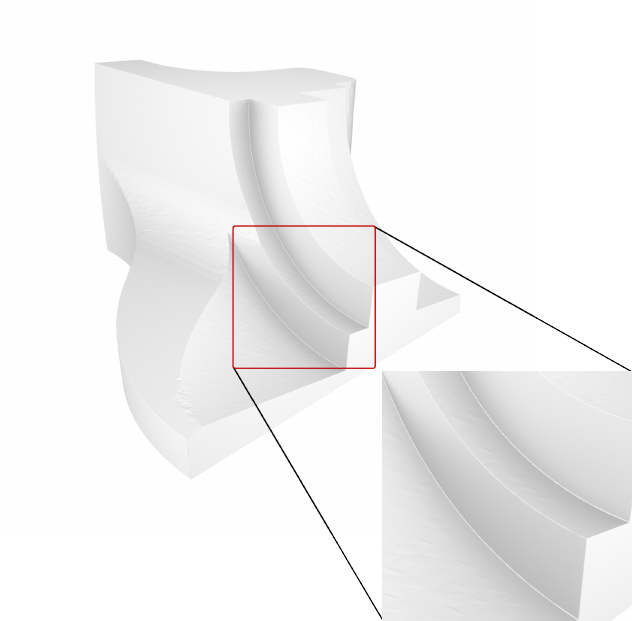}
	}
	\caption[]{Fan disk geometry reconstructed from the signed distance function~$\phi(\vec x)$ using bottom-up cut cell refinement and subsequent tetrahedralisation of cut quadrature cells.}   
	\label{fig:sharpFeatureCracks}
\end{figure*}

For nontrivial industrial CAD geometries, the bottom-up octree refinement of the cut cells leads to a significant improvement in the representation of the domain boundaries. The representation of the boundary of a fan disk with and without cut cell refinement are illustrated In Figure~\ref{fig:sharpFeatureCracks}.  In case of no refinement ($r_Q=0$), the surface has topological inconsistencies, and the sharp edges are poorly represented. In contrast, the boundary is faithfully  represented with refinement ($r_Q=3$). 

%-------------------------------------------------------------------------------
\subsection{Cut quadrature cell tetrahedralisation} 
\label{sec:cutTesselate}
%-------------------------------------------------------------------------------
%
We consider the decomposition of the quadrature cells~$\omega^Q_{i,k}$ that are cut by the domain boundary into tetrahedra. This is necessary to increase the accuracy of the evaluation of the finite element integrals. Whether~$\omega^Q_{i,k}$ is cut is determined by evaluating the signed distance~$\phi(\vec x)$ at its corners.  The part of a cut quadrature cell~$\omega^Q_{i,k}$ lying within the domain with~$\phi(\vec x) \ge 0$ is partitioned into several simplices such that 
\begin{equation} \label{eq:quadSimplices}
	 \omega^Q_{i,k} \Big \lvert_{\phi(\vec x) \ge 0}  	\approx \bigcup_l \tau_l \, ,
\end{equation}
see Figure~\ref{fig:sharpSubgridD}.
The finite element integrals are evaluated over the simplices~$\{  \tau_l \}$. Each cut quadrature cell~$\omega^Q_{i,k}$ has its own simplices. We did not make this dependence explicit in order not to clutter further the notation.  

\begin{figure*}[h!]
	\centering
	\includegraphics[scale=0.03]{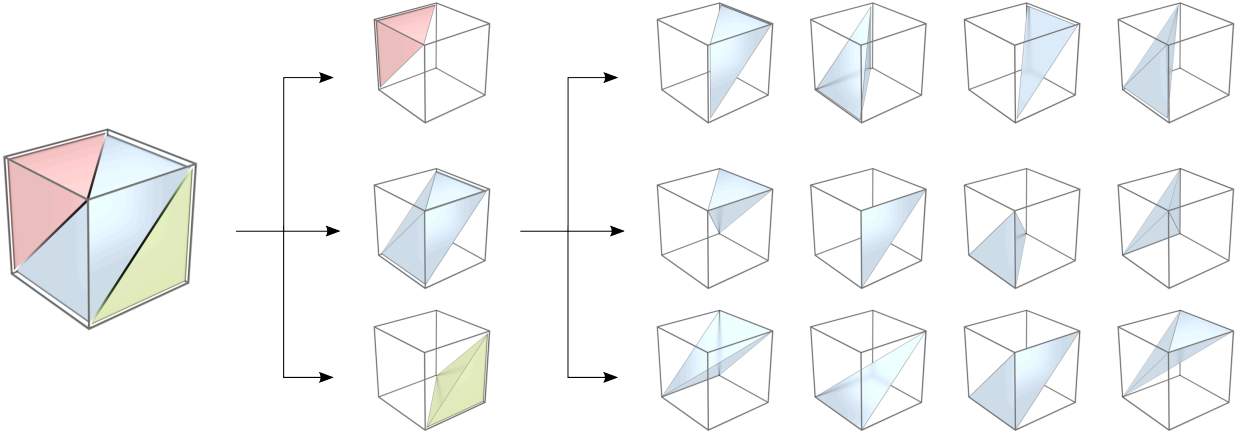}
	\caption{Splitting patterns for subdividing a cut quadtrature cell~$\omega_{i,k}^Q$ into six simplices prior to the applicationof the marching tetrahedra algorithm.}
	\label{fig:cubeSubdivision}
\end{figure*}

To obtain the simplices~$\{ \tau_l \}$, first, each cut cuboidal quadrature cell~$\omega^Q_{i,k}$ is split into six simplices and later the marching tetrahedra algorithm is used to decompose the part of the cell inside the domain~\cite{Gueziec:1995aa}.  The specific splitting pattern chosen in the first step is optimised so that the marching tethrahedra algorithm leads to as few simplices as possible. Specifically, as illustrated in Figure~\ref{fig:cubeSubdivision}, the cell is split into six simplices using  one of the three depicted splitting patterns~\cite{nielson1991asymptotic}. The domain boundary may cut all or some of the six resulting simplices. The tetrahedralisation of the parts of the cut simplices inside the domain is obtained with marching tetrahedra.  Depending on the sign of~$\phi(\vec x)$ at the corners of the cut simplices, there are three different possible splitting patterns in marching tetrahedra~\cite{Ruberg2012}. 
In the tetrahedralisation  of the cut simplices the intersection of their edges with the signed distance function~$\phi(\vec x)=0$ is required. To this end, we use a simple bisection algorithm. 
Note that the tetrahedralisation procedure also provides an approximation of the domain boundary using triangles within the FE cells, which are subsequently used for surface integrals.
%

%\begin{figure}[]
%	\centering
%	\subfloat[][Three nodes inside the domain \label{fig:marchingTetUniqueCasesA}]{
%		\includegraphics[scale=0.275]{./figs/immersed/marchingTetCasesA.pdf}
%	}	
%	\hspace{0.05\textwidth}	
%	\subfloat[][Two nodes inside the domain \label{fig:marchingTetUniqueCasesB}]{
%		\includegraphics[scale=0.275]{./figs/immersed/marchingTetCasesB.pdf}
%	}	
%	\hspace{0.05\textwidth}	
%	\subfloat[][One node inside the domain \label{fig:marchingTetUniqueCasesC}]{
%		\includegraphics[scale=0.275]{./figs/immersed/marchingTetCasesC.pdf}
%	}	
%	\caption{Three canonical splitting patterns of the marching tetrahedra algorithm. \label{fig:marchingTetUniqueCases}}
%\end{figure}

%
%-------------------------------------------------------------------------------
\subsection{Integration of cut quadrature cells}
\label{sec:feIntegration}
\begin{figure*}[]
	\centering
	\includegraphics[scale=1]{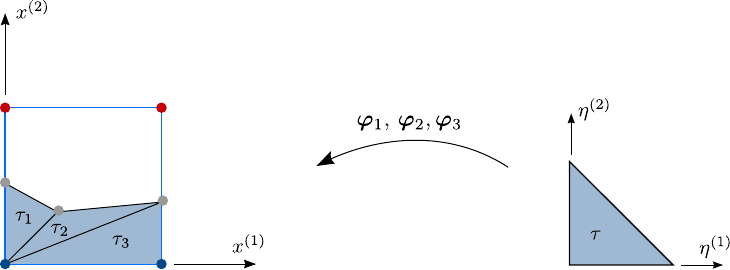}
	\caption{2D illustration of the integration of a quadrature cut cell~$\omega^Q_{i,k}$ corresponding to a finite element cut cell~$ \omega_{i} \in C^c$ . A reference simplex (right) is affinely mapped to the three simplices~$\tau_1$, $\tau_2$ and $\tau_3$ obtained with marching tetrahedra.}
	\label{fig:cutcellintegration}
\end{figure*}

%-------------------------------------------------------------------------------
%
We introduce now the evaluation of the finite element integrals over the cut quadrature cells~\mbox{$\omega^Q_{i,k} \subset \omega_i \in \set C^c$} in~\eqref{eq:quadCells} using the set of simplices~\mbox{$\{ \tau_l  \subset \omega^Q_{i,k}  \}$} in~\eqref{eq:quadSimplices}. Figure~\ref{fig:cutcellintegration} shows a typical setup with a cut quadrature cell and the associated integration simplices. The quadrature rules are given for a reference simplex~$\tau$ which is mapped to~$\tau_l$ using the affine mapping $\vec \varphi_l :  \vec \eta \in \tau \mapsto  \vec x \in \tau_l $.  Focusing, for instance, on the stiffness integral, its quadrature is given by 
\begin{equation} \label{eq:systemMatrix}
	\begin{aligned}	
		\int_{\omega^Q_{i,k}} \nabla u^h \cdot \nabla v^h \D \omega &\approx \sum_l \int_{\tau_l} \nabla u^h \cdot \nabla v^h \D \tau \\ &= 
		 \sum_l  \sum_g  \left(  \nabla u^h (\vec x_g )  \cdot \nabla  v^h (\vec x_g )  \right  ) \left  | \nabla_{\vec \eta} \vec \varphi_l  \right  |  w_g \, , 
	 \end{aligned}
\end{equation}	
where~$\vec x_g$ are the quadrature points mapped to the physical domain, $w_g$ are the weights, and~$  | \nabla_{\vec \eta} \vec \varphi_l   |$ is the absolute value of the determinant of the Jacobi matrix of the affine mapping. The solution~$u^h$ and the test function~$v^h$ are according to~\eqref{eq:interpolation} given in terms of the shape functions~$N_i(\vec x)$ of the finite element cell~$\omega_i$. Hence, the quadrature points~$\vec \eta_g$ must be mapped to the respective points~$\vec x_g$. This involves the  mapping~$\vec \varphi_l$ and the mapping implied by the bottom-up refinement of the finite element cell~$\omega_i$ into~$\omega^Q_{i,k}$.
% It is worth emphasising that both maps are affine so that the integrand in~\eqref{eq:systemMatrix} is polynomial and can be efficiently integrated. 
%
%

%-------------------------------------------------------------------------------
\subsection{Cut cell stabilisation}
\label{sec:stabilisation}
%-------------------------------------------------------------------------------
%
The cut finite element cells in~$\set C^c$  may have a small overlap with the physical domain~$\Omega$, and some of the associated basis functions may have a very small contribution to the system matrix~$\vec A$ in~\eqref{eq:discEquilibrium} resulting in an ill-conditioned matrix~\cite{dePrenter2017, de2023stability}. One approach to improving the conditioning is an elimination of their respective coefficients from the system matrix. Simply discarding some of the coefficients and basis functions would harm the finite element convergence rates.  Therefore, we  eliminate the critical coefficients by extrapolating the solution field from the nearby nodal coefficients of the cells  inside the domain. This approach is inspired by the extended B-splines by H\"ollig et al.~\cite{Hollig:2001aa}, which has also been applied in other immersed discretisation techniques~\cite{Ruberg2012,Ruberg2016, Xiao2019DG}.   

\begin{figure}[h!]
	\centering
	\includegraphics[scale=0.7]{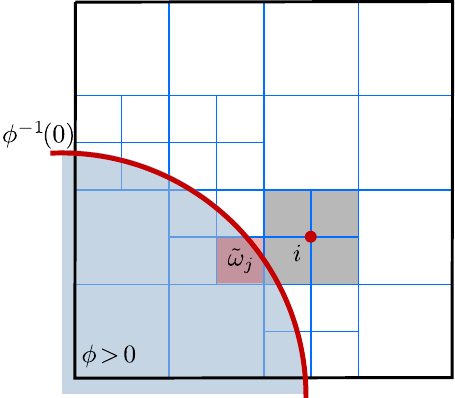}
	\caption[]{The support of the critical basis function $\supp N_i(\vec x)$ (in grey) corresponding to node~$i$ (in red) has a small overlap with the domain with~$\phi_h (\vec x) \ge 0$  (in blue). The respective coefficient~$u_i$ is expressed as the linear combination of the coefficients of the basis functions of the nodes of the active cell $\tilde{\omega}_j$. }   
	\label{fig:extBasisB}
\end{figure}

A basis function~$N_i (\vec x)$ is critical  when the intersection of its support with the physical domain~$\Omega$ lies below a threshold. Formally, this is expressed as 
\begin{equation}
\label{eq:criticality_criterion}
	\frac{\left | \supp  N_i(\vec x)   \cap \left  \{ \vec x \in \mathbb R^3 \mid \phi (\vec x) \ge 0 \right  \} \right  |  }{ | \supp  N_i(\vec x)  | }   < \epsilon \,,  
\end{equation}
where~$| \cdot |$ is the volume of the respective set, and the threshold is chosen as~$\epsilon = 1/8$ in our numerical computations, see Figure~\ref{fig:extBasisB}.
To determine the numerator and denominator of (\ref{eq:criticality_criterion}) for each basis function $N_i(\vec x)$, we perform an assembly-like procedure, in which each FE cell contributes its active and total volume to the corresponding basis functions.

The coefficients~$u_i$ of a critical basis function~$N_i(\vec x)$ are extrapolated from the nodal coefficients of a cell~$\widetilde{\omega}_j$ which contains only non-critical basis functions. To identify $\widetilde{\omega}_j $ we first collect all the active cells in the two neighbourhood of the node~$\vec x_i$ corresponding to~$N_i(\vec x)$.
%\ADDED{That is, we collect elements contributing to the basis function and their neighbours over faces, edges and vertices.}
Subsequently, we select from the candidate cells the cell~$\widetilde{\omega}_j$ with the centroid closest to node~$i$. Denoting the set of nodal indices of ~$\widetilde{\omega}_j$  with~$\mathbb J$  the critical coefficient~$u_i$ is obtained by simply evaluating the basis functions~$N_j(\vec x)$  with $j \in \mathbb J$ at~$\vec x_i$, that is, 
\begin{equation}
\label{eq:extrapolation_coeffs}
    u_i = \sum_{j \in \mathbb J } N_j (\vec x_i) u_{j} \, . 
\end{equation}
After extrapolating the coefficients of all the critical basis functions present in the grid, the original nodal coefficients can be expressed abstractly as 
\begin{equation}
\label{eq:change_of_basis}
	\vec u = \vec E \, \widetilde{ \vec u} \, .
\end{equation}
Introducing this relation in the discrete system of equations~\eqref{eq:discEquilibrium} we obtain
\begin{equation}
\label{eq:transformed_system}
	\widetilde{ \vec A}\, \widetilde{ \vec u} = \widetilde{ \vec f} \,, 
\end{equation}
where $\widetilde{ \vec A} = \vec E^\trans \vec A \vec E$ and $\widetilde{ \vec f} = \vec E^\trans  \vec f$. The matrix~$\widetilde{ \vec A}$ is well conditioned so that~\eqref{eq:transformed_system} can be robustly solved.

%\ADDED{In some cases, the active cells in the two neighbourhood of node~$\vec x_i$ may belong to different partitions and computing cores. We select for extrapolation a cell~$\widetilde{\omega}_j $ which lies in the same partition like~$\vec x_i$. 
%When no suitable candidate exists within the same partition, any cell that contains the critical basis function~$N_i(\vec x)$ is excluded from the assembly. This could happen in a rare situation where the critical basis function lies in the proximity of the interface between subdomains.
%d}

%
%

	%--------------------------------------------------------------------------------
\section{Domain partitioning and parallel solution}
\label{sec:parallel}
%--------------------------------------------------------------------------------

%--------------------------------------------------------------------------------
\subsection{Finite element grid partitioning}
\label{sec:partitioning}
%--------------------------------------------------------------------------------

In this section, we describe the decomposition of the FE grid into a set of subdomains assigned to individual processors of a distributed parallel computer. As mentioned in Section~\ref{sec:femGrid}, the FE grid is generated by refining the base grid through selective octree refinements. We can either perform refinements towards the geometric boundary, i.e., by identifying the cut cells $\set C^c$, or towards the solution features. The latter is more demanding since the problem has to be solved repeatedly, or otherwise requires a-priori estimates of the solution.

We consider a division of the FE grid in $\Omega_\Box $ into $N_S$ nonoverlapping subdomains, $\Omega_{\Box i},\ i = 1,\dots,N_S$, which are assigned to individual processor cores and further processed in parallel. This is achieved by first generating the Z-curve of the FE grid, and subsequently subdividing the Z-curve into partitions of approximately equal length. It is important to note that each FE cell has predicate active, cut, or inactive, where only the first two kinds contribute significantly to the computation. Considering these characteristics in the partitioning leads to a more balanced load distribution in each processor. Specific to the processing of the cut-cell identification algorithm, we restrict the face-neighbour search within each subdomain.

\begin{figure*}[h!]
	\centering
	\subfloat[][\label{fig:withoutWeight}]{
		\includegraphics[scale=0.16]{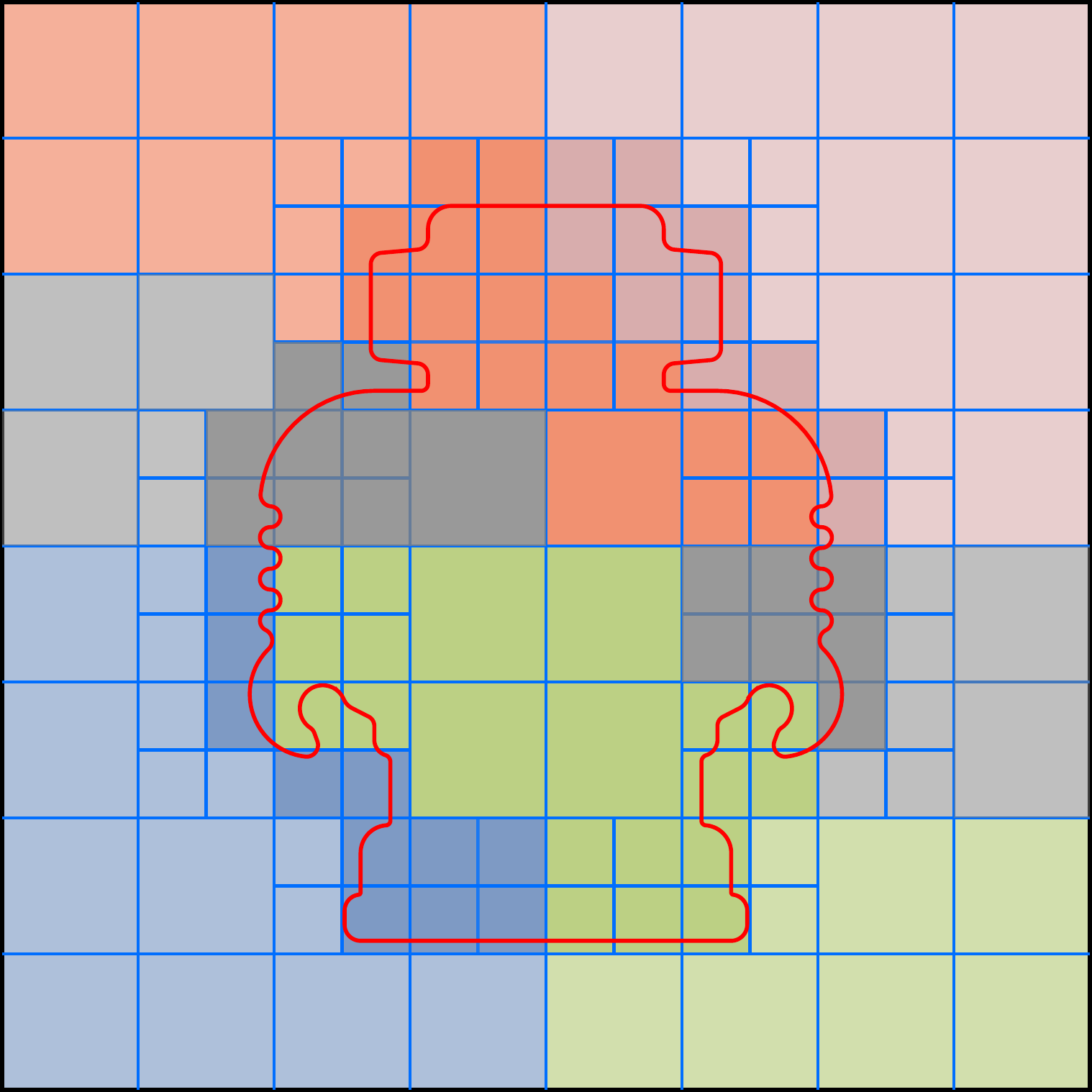}
	}
	\hspace{0.1\textwidth}	
	\subfloat[][\label{fig:withWeight}]{
		\includegraphics[scale=0.16]{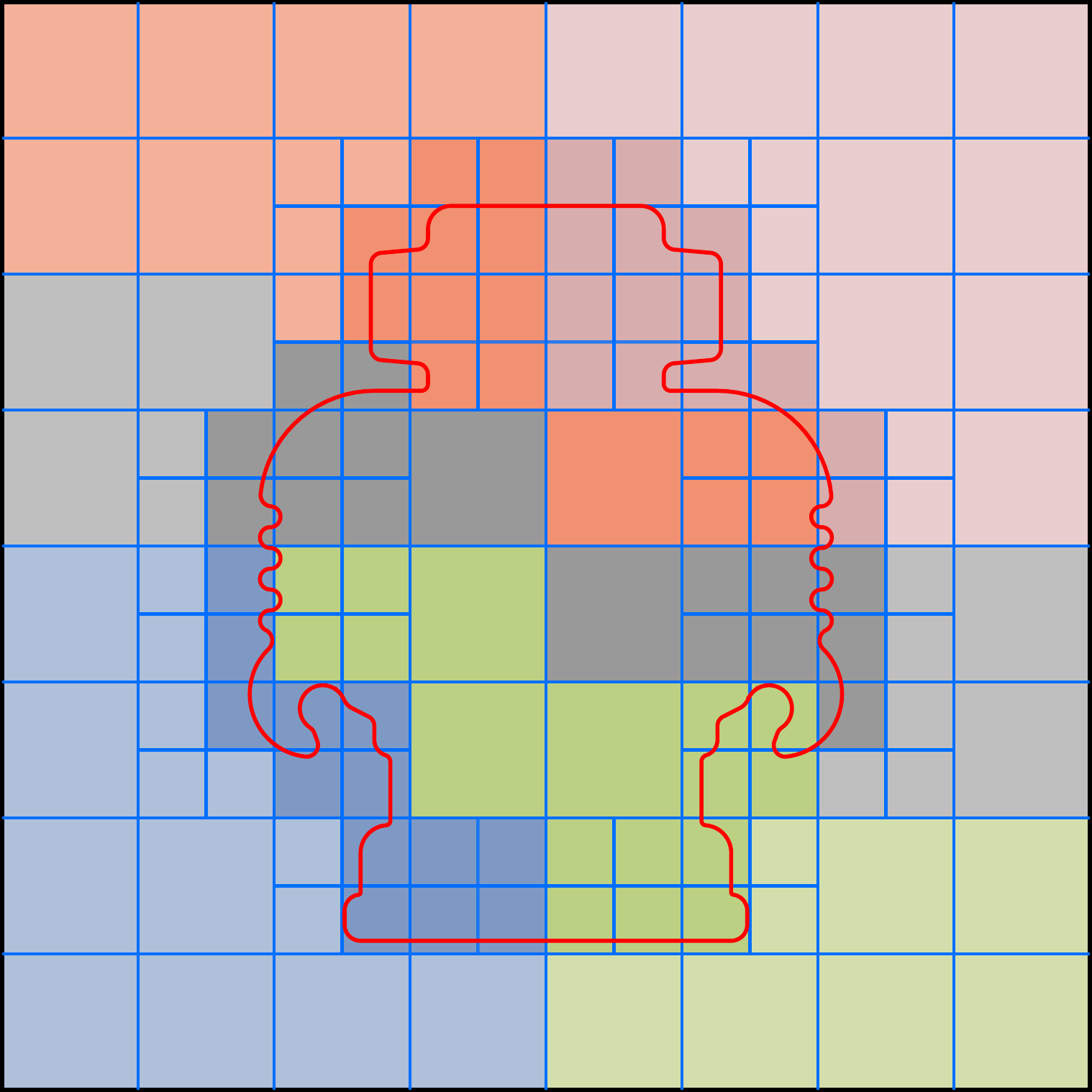}
	}
	\caption[Decomposition of the domain $\Omega_{\square}$]{Domain decomposition of $\Omega_{\square}$. Darker shades indicate that the cell is either active or cut, i.e. it belongs to $\set C^d$. (a) Domain decomposition based solely on equal partition of the Z-curve. The total number of cells of different colours are roughly equal but the active and cut cells differ for each colour. (b) A more balanced domain decomposition is achieved by applying larger weights for active and cut cells. }   
	\label{fig:domPartition}
\end{figure*}

In particular, we consider applying different weights on FE grid cells in the paritioning of the Z-curve. Let us denote the union of the active cells and cut cells as $\set C^d = \set C^a \cup \set C^c$. We apply for the set $\set C^d$ a larger weight of 100 than those for inactive cells $\set C^i$ of weight 1. Then, the division into subdomains is inherited from the division of $\Omega_\Box $ as intersections of $\set C^d$ with subdomains $\Omega_{\Box i}$, i.e. $\Omega_i = \set C^d \cap \Omega_{\Box i},\ i = 1,\dots,N_S$. The weighted subdivision yields a decomposition of the domain with a balanced number of cut and active cells per processor, as illustrated in Figure~\ref{fig:domPartition}. 

Depending on the computing memory requirements of the application, the geometry grid can be replicated on each partition or partitioned using the level set fracturing functionality of openVDB according to the bounding box of each partition~\cite{openvdbwebpage}. The quadrature grid is generated in the cut cells belonging to each partition according to Section~\ref{sec:bottomUpCut} by sampling the level set values. 

%--------------------------------------------------------------------------------
\subsection{Hanging nodes}
\label{sec:hanging}
%--------------------------------------------------------------------------------

\begin{figure}[htp]
	\centering
	\includegraphics[width=0.55\linewidth]{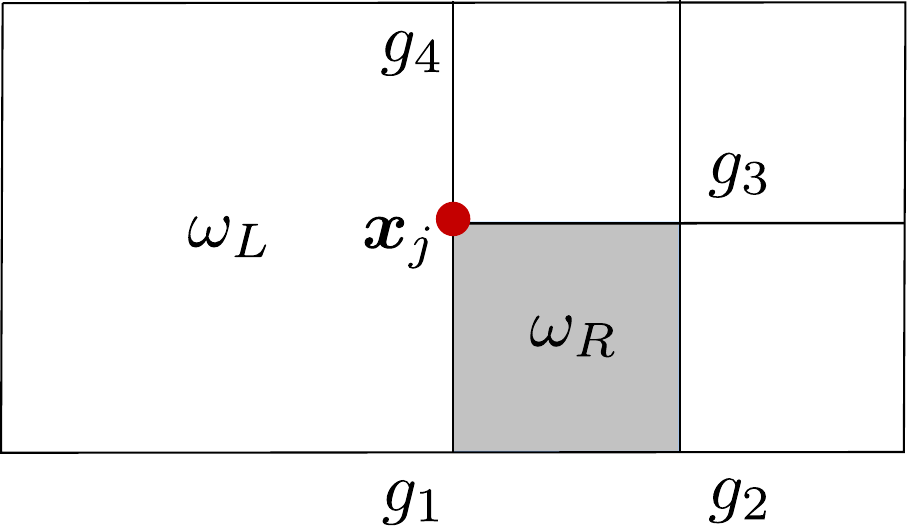}
 	\caption{Example of a hanging node. Large constraining element $\omega_L$ determines the value in the hanging node (red dot) of the constrained element $\omega_R$ (gray). Hanging node does not have a degree of freedom in the global system, its value is dictated by values of degrees of freedom $g_1$ and $g_4$. Algorithmically the constraint can be realized by applying the change-of-basis matrix $\vec T_R$ to the element matrix of $\omega_R$ and assembling the transformed local matrix to degrees of freedom $g_1$, $g_2$, $g_3$, and $g_4$.}   
	\label{fig:hanging}
\end{figure}

Octree refinement in constructing the FE grid lead to the presence of hanging nodes.
These are nodes created at the faces and edges between two or more elements by refining only one of them, see Figure~\ref{fig:hanging}. For preserving the continuity of the FE solution, we require that the solution at the hanging node is determined by the values in the nodes of the larger element.

Consider the constraining element $\omega_L$ and the constrained element $\omega_R$, as exemplified in Figure~\ref{fig:hanging}. We express the solution at the hanging node at $\vec x_j$ as 
\begin{equation}
\label{eq:hangingConstraints}
     u_h(\vec x_j) = \sum_{k = 1}^{n_k} H_{j,k} u_{k} \, ,
\end{equation}
where $H_{j,k}$ is the value of the $k-$th shape function associated with the element $\omega_L$ at $\vec x_j$, i.e., $H_{j,k} = N_k(\vec x_j)$. Given that the only nonzero shape functions are those associated to the nodes at the common edge, the solution at the hanging node is constrained by the nodal values of $\omega_L$ at the edge, i.e., at the degrees of freedom $g_1$ and $g_4$. In our example with the bilinear shape functions, (\ref{eq:hangingConstraints}) simplifies to $u_h = \frac{1}{2}u_1 + \frac{1}{2}u_4$, where $u_1$ and $u_4$ are the values of the degrees of freedom $g_1$ and $g_4$, respectively. With respect to the implementation, the procedure for eliminating hanging nodes is equivalent to the one for eliminating the critical nodes through extrapolation in~\eqref{eq:extrapolation_coeffs}.

As a final remark, we note that it is possible for a hanging node to be constrained by a node which is subject to extrapolation. In the extreme case, the node from which we extrapolate can also be a hanging node and constrained by another regular node. This situation is in fact handled naturally by nesting the assembly lists related to the hanging nodes and those related to the extrapolation. At the end of this nesting, the element contributes its local matrix to regular degrees of freedom potentially through a relatively long assembly list.
A more formal description of the potential chaining of constraints due to hanging nodes and the extrapolation has been recently given in~\cite{badia2021aggregated}.

%--------------------------------------------------------------------------------
\subsection{Iterative substructuring}
\label{sec:Substructuring}
%--------------------------------------------------------------------------------

We employ iterative substructuring method to solve the linear system arising from the immersed finite element system~\eqref{eq:transformed_system}. We consider that the global stiffness matrix $\widetilde{\vec A}$ and the right-hand side vector $\widetilde{\vec f}$ appearing in~\eqref{eq:transformed_system} can be assembled from the local contributions from each subdomain $\Omega_{\Box i}$, that is, $ \widetilde{\vec A} = \sum_{i=1}^{N_S} \vec R_i^\trans \vec A_{i}  \vec R_i$ and $\sum_{i=1}^{N_S} \widetilde{\vec f}  = \vec R_i^\trans \vec f_{i}$. The Boolean restriction matrix $\vec R_{i}$ containing a single~$1$ in each column selects the local from the global degrees of freedom. 

For each local subdomain, we separate the degrees of freedom belonging to interior $\vec u_i^{I}$ and the interface $\vec u_i^{\Sigma}$, which leads to a 2$\times$2 blocking of the local linear system,
\begin{equation}
	\label{eq:subdomainBlocking}
	\left[
	\begin{array}[c]{cc}
		\vec A_{i}^{II}       & \vec A_{i}^{I\Sigma} \\
		\vec A_{i}^{\Sigma I} & \vec A_{i}^{\Sigma\Sigma}  
	\end{array}
	\right]
	\left[
	\begin{array}[c]{cc}
		\vec u_{i}^{I} \\
		\vec u_{i}^{\Sigma}  
	\end{array}
	\right]
	=
	\left[
	\begin{array}[c]{cc}
		\vec f_{i}^{I} \\
		\vec f_{i}^{\Sigma}  
	\end{array}
	\right] \, .
\end{equation}
Here the local interface degree of freedom $\vec u_{i}^{\Sigma} $ can be assembled into a global set of unknowns at the interface of all subdomains $\Sigma$ utilising an interface restriction matrix $\vec R^{\Sigma}_i:\Sigma \rightarrow \Sigma_{i}$, where $\vec u^{\Sigma} = \sum_{i=1}^{N_S} \vec R_i^{\Sigma^ \trans} \vec u_{i}^{\Sigma} $.
Interface $\Sigma$ is formed by degrees of freedom shared by several subdomains.

The iterative substructuring seeks the solution of the global interface problem 
\begin{equation}
	\label{eq:Sug}
	\vec S \vec u^{\Sigma} = \vec h \, , 
\end{equation}
using, for instance, the preconditioned conjugate gradient (PCG) method. The global Schur complement matrix $\vec S$ comprises of the subdomain contribution
\begin{equation}
    \vec S = \sum_{i=1}^{N_S} \vec R_i^{\Sigma^\trans} \vec S_{i} \vec R^{\Sigma}_i \, ,
\end{equation}
%$\vec S = \sum_{i=1}^{N_S} \vec R_i^{\Sigma \trans} \vec S_{i} \vec R^{\Sigma}_i$, 
where the local Schur complement with respect to $\Sigma _i$ is defined as
\begin{equation}
	\label{eq:S-assembly}
	\vec S_i = \vec A_{i}^{\Sigma\Sigma} - \vec A_{i}^{\Sigma I} \left(\vec A_{i}^{II} \right)^{-1} \vec A_{i}^{I\Sigma} \, . 
\end{equation}
Similarly, the right-hand side is assembled from the subdomains 
\begin{equation}
\vec h  = \sum_{i=1}^{N_S} \vec R_i^{\Sigma^\trans} \vec h_{i} \, , 
\end{equation}
where 
\begin{equation}
	\label{eq:h-assembly}
	\vec h_i = \vec f_i^{\Sigma} - \vec A_{i}^{\Sigma I} \left(\vec A_{i}^{II}\right)^{-1} \vec f_{i}^{I}.
\end{equation}
Once we know the local solution at the interface $\vec u_i^{\Sigma}$, the solution in the interior of each subdomain $\vec u_{i}^{I}$ is recovered from the first row of \eqref{eq:subdomainBlocking}. Note that neither the global matrix $\vec S$ nor the local matrices $\vec S_i$ are explicitly constructed in the iterative substructuring. Only multiplications of vectors with $\vec S_i$ are needed at each PCG iteration. 

%--------------------------------------------------------------------------------
\subsection{BDDC preconditioner}
\label{sec:BDDCpreconditioner}
%--------------------------------------------------------------------------------

Next, we briefly describe the balancing domain decomposition based on constraints (BDDC) preconditioner in solving the interface problem (\ref{eq:Sug}).
An action of the BDDC preconditioner $\vec M_{BDDC}^{-1}$ produces a preconditioned residual $\vec z^{\Sigma}$ from the residual in the $k$-th iteration 
$\vec r^{\Sigma} = \vec S \vec u^{\Sigma}_{(k)} - \vec h$
by implicitly solving the system $\vec M_{BDDC} \, \vec z^{\Sigma} = \vec r^{\Sigma}$.
Specifically, BDDC considers a set of coarse degrees of freedom which will be continuous across subdomains such that the preconditioner is invertible yet inexpensive to invert.
In this work we consider degrees of freedom at selected interface nodes (corners) and arithmetic averages across subdomain faces and edges as the coarse degrees of freedom.
This gives rise to a global coarse problem with the unknowns $\vec u_{C}$ and local subdomain problems with independent degrees of freedom $\vec u_{i}$ which are parallelisable.

BDDC gives an approximate solution which combines the global coarse and local subdomain components,
i.e.,
\begin{equation} 
	\label{eq:averaging_subdomain_solves}
	\vec z^{\Sigma} = \sum_{i=1}^{N_S} \vec R_{i}^{\Sigma^\trans} \vec W_{i} \vec R_{Bi} \left( \vec u_{i} + \vec \Phi_{i} \vec R_{Ci} \vec u_{C} \right).
\end{equation} 
In particular, $\vec u_{C}$ and $\vec u_{i}$ are obtained in each iteration by solving
\begin{align} 
	\label{eq:coarse_problem}
	&\vec S_{C} \vec u_{C} = \sum_{i = 1}^{N_S} \vec R_{Ci}^{\trans} \vec \Phi_{i}^{\trans} \vec R_{Bi}^\trans \vec W_i \vec R_{i}^{\Sigma} \vec r^{\Sigma}, \\
	\label{eq:localCorr}
	&\left[
	\begin{array}[c]{cc}
		\vec A_{i} & \vec C_{i}^{\trans} \\
		\vec C_{i} &  \vec 0  
	\end{array}
	\right]  \left[
	\begin{array}[c]{c}
		\vec u_{i}\\
		\vec \mu_{i}
	\end{array}
	\right] = \left[
	\begin{array}[c]{c}
		\vec R_{Bi}^\trans \vec W_i \vec R_{i}^{\Sigma} \vec r^{\Sigma}\\
		\vec 0
	\end{array}
	\right]\, , \, i=1,\dots,N_{S},
\end{align} 
where $\vec S_{C}$ is the stiffness matrix of the global coarse problem, matrix $\vec A_{i}$ is assembled from elements in the $i$-th subdomain, and $\vec C_{i}$ is a constraint matrix enforcing zero values of the local coarse degrees of freedom. The diagonal matrix $\vec W _i$ applies weights to satisfy the partition of unity, and it corresponds to a simple arithmetic averaging in this work. The Boolean restriction matrix $\vec R_{Bi}$ selects the local interface unknowns from those at the whole subdomain, the columns of $\vec \Phi_{i}$ contain the local coarse basis functions,
and $\vec R_{Ci}$ is the restriction matrix of the global vector of coarse unknowns to those present at the $i$-th subdomain. Our implementation of the BDDC preconditioner is detailed in~\cite{Kus-2017-CND}. For the interested readers, we refer to~\cite{Dohrmann-2003-PSC, Fragakis2003, Toselli-2005-DDM} for a thorough description of the BDDC method and its multilevel variants~\cite{Tu-2007-TBT3D, Mandel-2008-MMB, Badia-2018-RSD}.

\begin{figure}[htp]
	\centering
	\includegraphics[width=0.55\linewidth]{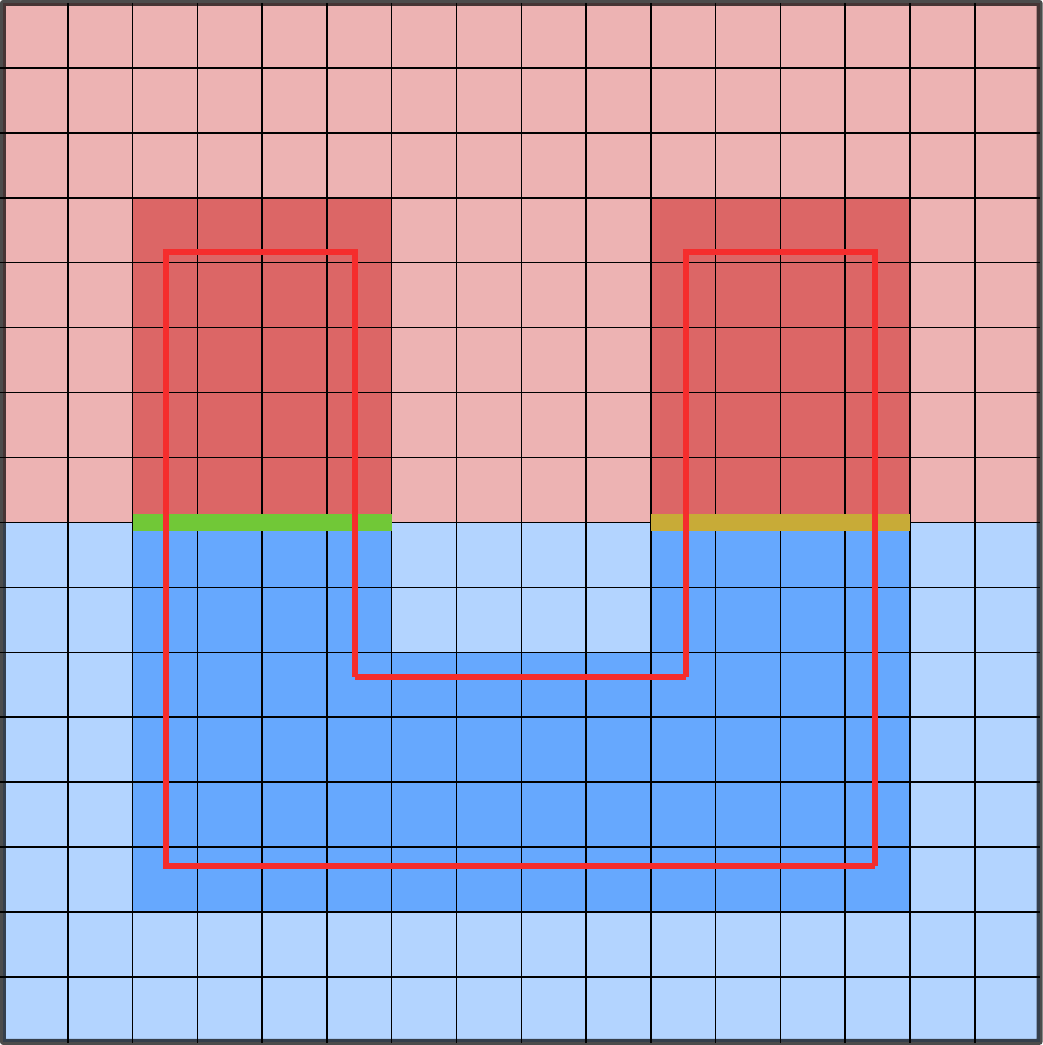}
 	\caption{Example of fragmenting of subdomains due to the extraction of the active and cut cells $\set C^d$. By analysis of the dual graph of each subdomain, the two components of the red subdomain are detected, and the two faces between red and blue subdomains are identified.
	\label{fig:fragmenting}
}
\end{figure}

In the context of the domain decomposition of our FE grid, we might encounter the issue of subdomain fragmentation. One typical cause is the partitioning of the Z-curve which does not guarantee connected subdomains, see for example the gray subdomain showcased in Figure~\ref{fig:withoutWeight}. Another possible cause of the fragmenting is due to the subdomain extraction $\Omega_i = \set C^d \cap \Omega_{\Box i}$. A simple case of this effect is shown in Figure~\ref{fig:fragmenting}.

As proposed in~\cite{Kus-2017-CND}, we remedy this problem by analyzing the \emph{dual graph} of the FE grid of each subdomain followed by generation of inter-subdomain faces and corresponding constraints independently for each component.
Within the dual graph, FE grid cells correspond to graph vertices, and a graph edge is introduced between two vertices whenever the corresponding cells share at least four degrees of freedom. Consequently, the local saddle point problem required for obtaining the local basis functions $\vec \Phi_i$ will become solvable. In our implementation, we rely on the open-source \textsl{BDDCML} solver~\cite{Sousedik-2013-AMB} for the implementation of the multilevel BDDC method, and the solver is equipped with a component analysis based on the dual graph of each subdomain grid.

	%--------------------------------------------------------------------------------
\section{Examples \label{sec:examples}}
%--------------------------------------------------------------------------------
%\todo[inline]{Pavel: Maybe we could start with a quick recap of some of the parameters we will use in the following, at least $h_g$, $h_q$, $h_f$...}

We introduce several examples of increasing complexity to demonstrate the convergence of the proposed immersed finite element scheme and its robustness for complex 3D CAD geometries. We focus the analysis on the Poisson problem, e.g., modelling heat transfer, where the solution is a scalar field. Throughout this section, we emphasise the use of three grids characterised by their resolution, namely $h_g$ for geometry grid, $h_f$ for FE grid and $h_q$ for quadrature grid.

%--------------------------------------------------------------------------------
\subsection{Interplay between geometry, quadrature and FE grid sizes \label{sec:convergence_sphere}}
\begin{figure}[htp]
	\centering
	\includegraphics[width=0.25\textwidth]{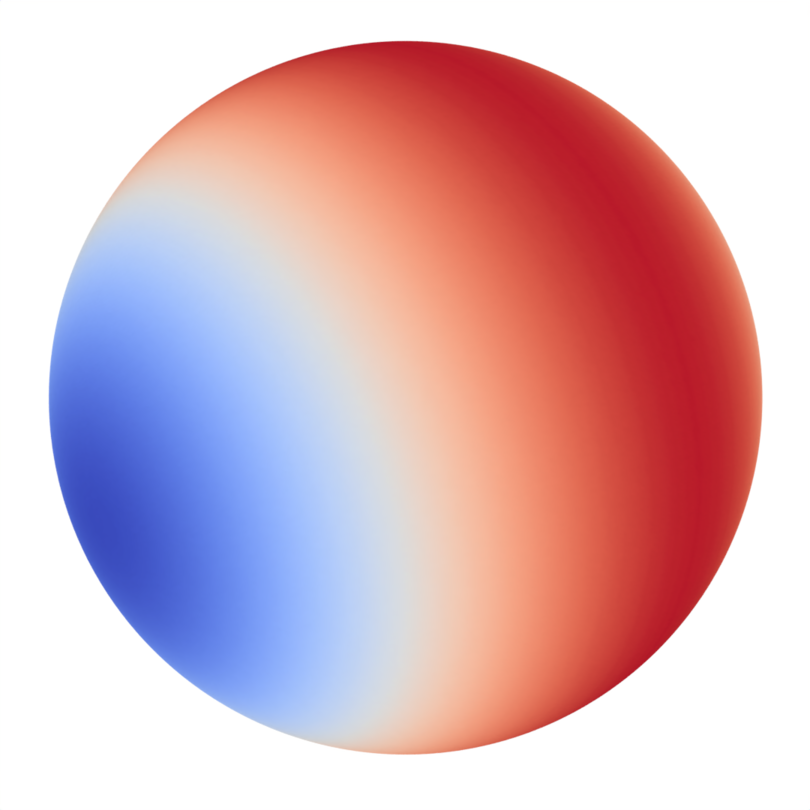} 
	\caption{Solution of the Poisson problem in a sphere, $u = \cos(x_3)$.}
	\label{fig:sphere_solution}
\end{figure}

%--------------------------------------------------------------------------------
%
%\todo[inline]{$L_2 $ and $H_1$ norm convergence for linear and higher order basis function}
%\todo[inline]{Can choose a rather standard geometry - fandisk or bunny}
%\todo[inline]{Comparison of cost between adaptive and non-adaptive}
%\begin{equation}
%	\label{eq:cosine}
%	u = \cos(x_3) \, ,
%\end{equation}
As a first example we consider the Poisson-Dirichlet problem on a unit sphere with a prescribed solution $u = \cos(x_3)$ as illustrated in Figure~\ref{fig:sphere_solution}. 

\begin{figure*}[]
	\centering
 	\subfloat[ ][$h_q=0.6$]{
	   \includegraphics[width=0.25\textwidth]{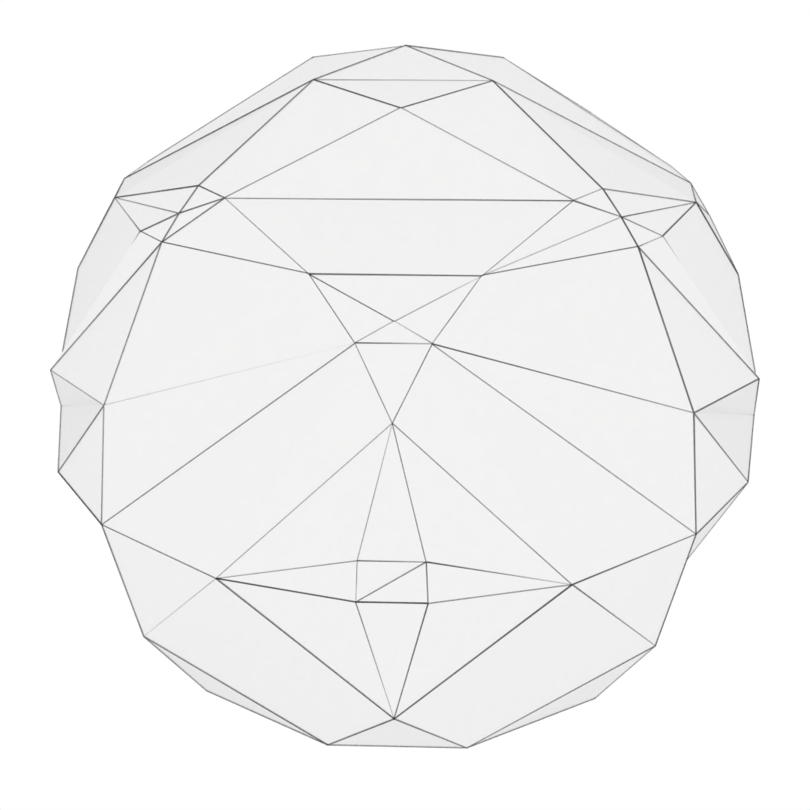} }\hspace{0.02\textwidth}
    \subfloat[][$h_q=0.15$]{
        \includegraphics[width=0.25\textwidth]{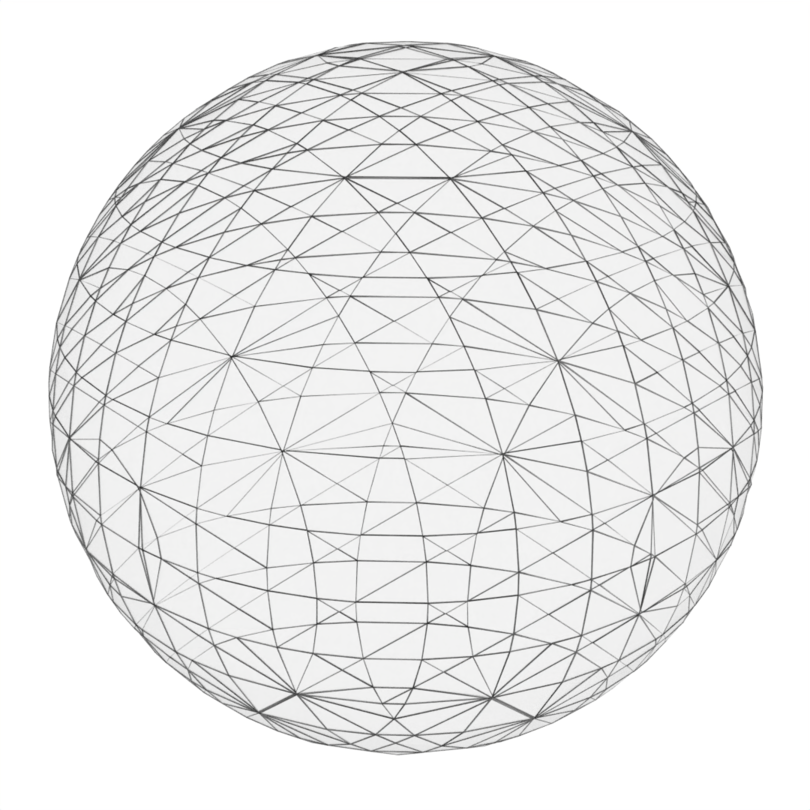} }\hspace{0.02\textwidth}
    \subfloat[][$h_q=0.035$]{    
	\includegraphics[width=0.25\textwidth]{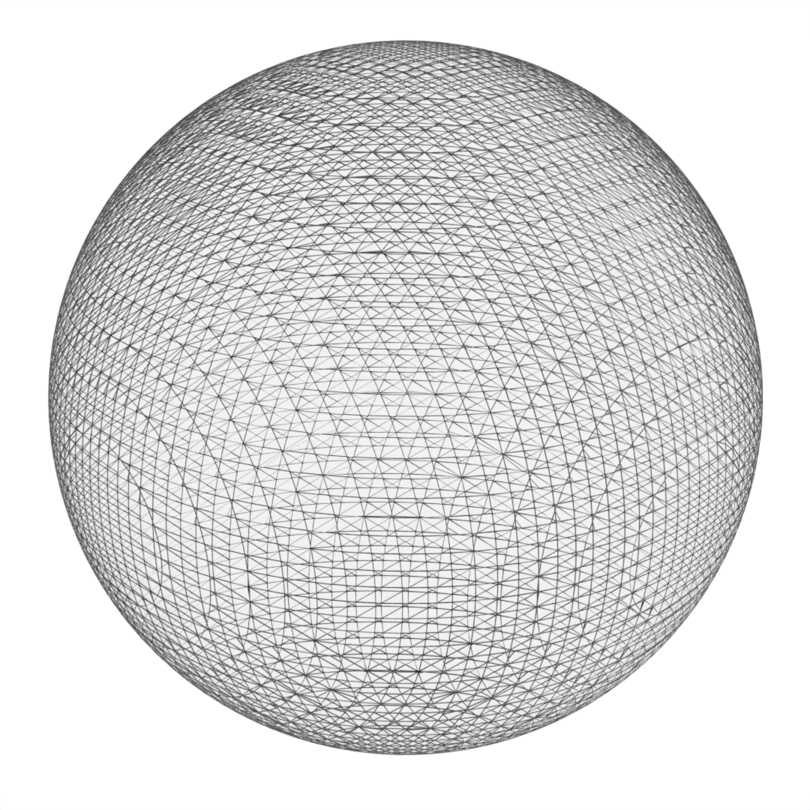} 
    }
	\caption{Sphere reconstructed from the quadrature grid using the marching tetrahedra algorithm for a fixed geometry grid resolution~$h_g=1/128$ and different quadrature grid resolutions~$h_q$.}
	\label{fig:sphere_reconst}
\end{figure*}

\begin{figure*}[htp]
	\centering
	\subfloat[][Volume]{
		\includegraphics[width=0.47\textwidth]{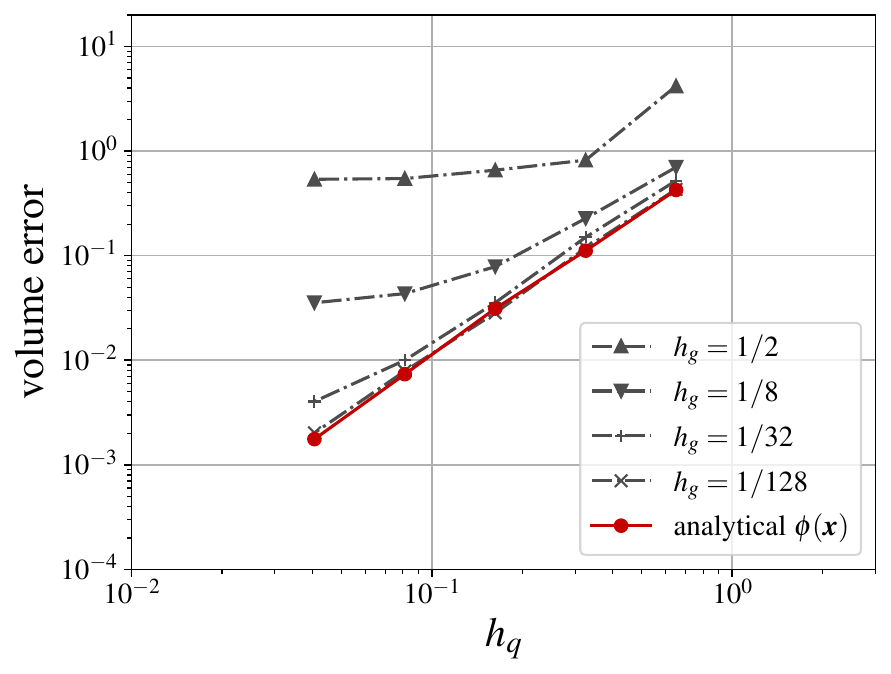} 
	}
	\hspace{0.02\textwidth}%
	\subfloat[][Surface]{
		\includegraphics[width=0.47\textwidth]{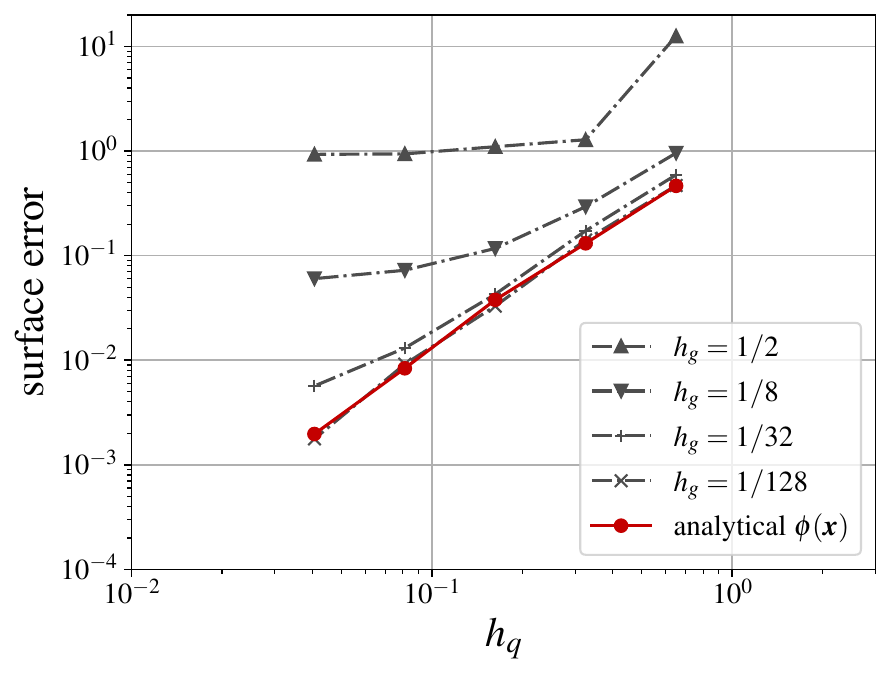}	
	}
	\caption{Dependence of the error in volume and surface computation for a fixed quadrature grid resolution~$h_q=0.035$ and different geometry grid resolutions~$h_g$. The red lines indicate the surface and volume errors when analytical signed distance function is used.}
	\label{fig:volume_surface}
\end{figure*}

In this example, we first study the effect of the geometry grid resolution $h_g$ and quadrature grid resolution $h_q$ on the volume and surface integration. The sphere geometry is described using an analytical signed distance function $\phi(\vec x) = 1 - (x_1^2 + x_2^2 + x_3^2)^{1/2}$. 
To start with, we sample the signed distance function over the nodes of the geometry grid with a resolution $h_g = \{1/2, \, 1/8, \, 1/32, \, 1/128 \}$.
%As can be seen in Figure~\ref{fig:sphere_levelset}, finer resolution $h_g$ yields a more faithful representation of the exact geometry.
For each value of $h_g$, we use quadrature grid with resolution $h_q = \{0.6, \, 0.3, \, 0.15, \, 0.07, \, 0.035 \}$ to linearly reconstruct the volume and the surface by means of marching tetrahedra. Figure~\ref{fig:sphere_reconst} indicates the reconstructed sphere for various $h_q$ for the geometry grid resolution $h_g = 1/128$. We compare the volume and surface area of the reconstructed representation with the analytical value in Figure~\ref{fig:volume_surface}. For each $h_g$ it is evident that the errors reach plateau when $h_q < h_g$, indicating that the geometric error is bounded by the geometry grid resolution $h_g$. 

%The first experiment aims at the relation between the geometry grid resolution $h_g$ and the quadrature grid resolution $h_q$.
%It is tested by using an STL file with the sphere with a very high resolution and representing it on the grid with changing $h_g$.
%A sequence\ of level set accuracy for different values of $h_g$ is shown in Figure~\ref{fig:sphere_levelset}.\todo{Pavel: When the meshes are plotted in that figure, what was the $h_q$? Or what is the way it was plotted?}
%The accuracy of the representation is evaluated by the error in volume and surface compared to the analytic value.
%These plots are shown for different $h_q$ in Figure~\ref{fig:volume_surface}.
%A sequence of quadrature meshes\todo[inline]{Pavel: did we emphasise that we use the quadrature meshes for visualisation as well?} obtained by the marching tetrahedra algorithm for different values of $h_q$ is shown in Figure~\ref{fig:sphere_reconst}.\todo{Pavel: What was the $h_g$ there? Or was it analytic?}

%\begin{figure*}[]
%	\centering
%	\includegraphics[width=0.25\textwidth]{./figs/examples/sphere_lv_05} \hspace{0.02\textwidth}
%	\includegraphics[width=0.25\textwidth]{./figs/examples/sphere_lv_0125}  \hspace{0.02\textwidth}
%	\includegraphics[width=0.25\textwidth]{./figs/examples/sphere_lv_003}
%	\caption{Representation of $\Omega_g$ and $\Gamma_g$ for a sphere and different values of $h_g$: $1/2$ (left), $1/8$ (centre), and $1/32$ (right).}
%	\label{fig:sphere_levelset}
%gmai\end{figure*}

\begin{figure*}[htp]
	\centering
	\subfloat[][Relative $L_2$ norm error]{
		\includegraphics[width=0.47\textwidth]{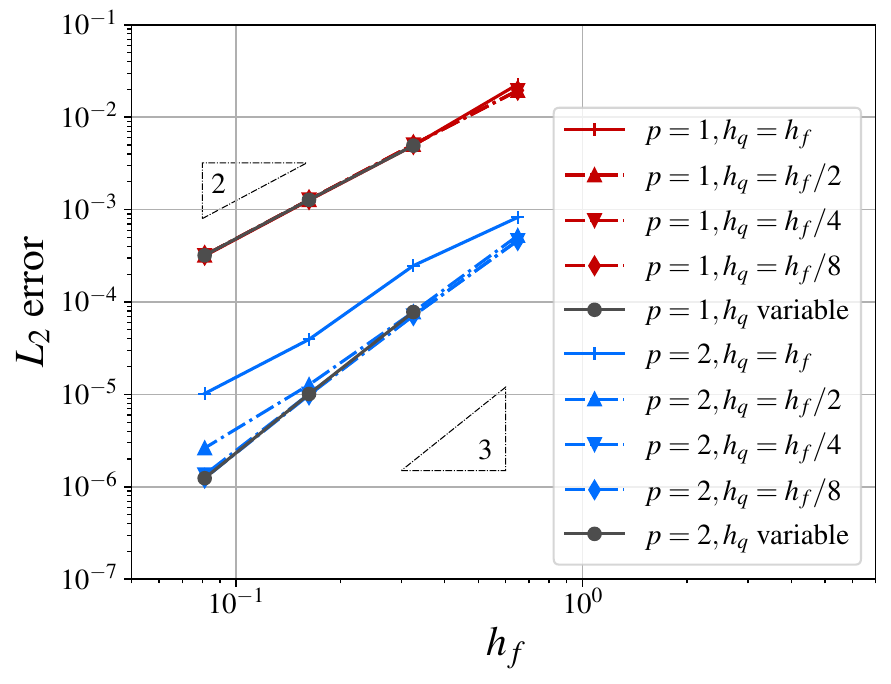} 
	}
	\hspace{0.02\textwidth}%
	\subfloat[][Relative $H^1$-seminorm error]{
		\includegraphics[width=0.47\textwidth]{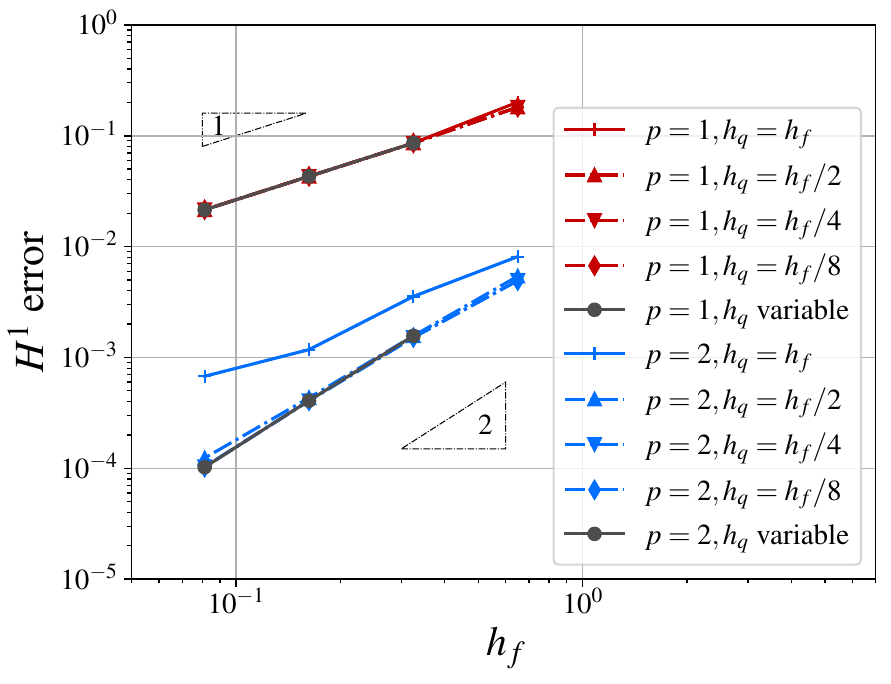}	
	}
	\caption{Dependence of the relative finite element errors for linear ($p=1$) and quadratic basis ($p=2$) functions and for different finite element and quadrature grid resolutions $h_f$ and $h_q$.}
	\label{fig:sphere_convergence_h}
\end{figure*}

Next, we verify the numerical convergence of the relative error for the Poisson-Dirichlet problem both in $L_2$-norm, i.e. \mbox{$\|u-u_h\|_{L_2}/\|u_h\|_{L_2}$}, and $H^1$-seminorm, i.e. \mbox{$|u-u_h|_{H^1}/|u_h|_{H^1}$}.  We take a fine STL mesh of the sphere to obtain the signed distance function $\phi(\vec x)$. The resolution of the geometry grid is twice smaller than the finest resolution of the quadrature grid, i.e., $h_g = h_q / 2$, where the finest quadrature grid resolution is $h_q = h_f / 8$. The convergence is established against the resolution of the uniform FE grid $h_f$. Note that FE grid is used for the field discretisation using Lagrange basis functions. For integration we use a quadrature grid of size \mbox{$h_q = \{h_f, \, h_f / 2, \, h_f / 4, \, h_f / 8  \}$}. As an example, $h_q = h_f / 4$ indicates that quadrature grid is obtained by twice subrefining the cut FE cells. We also include the case $h_q \equiv h_f$, i.e., $h_q$ is refined whenever $h_f$ is refined. Figure~\ref{fig:sphere_convergence_h} shows the convergence in $L_2$-norm and $H^1$-seminorm for both linear ($p = 1$) and quadratic ($p = 2$) basis functions. It is evident from Figure~\ref{fig:sphere_convergence_h} that for the linear case optimal convergence rate is achieved for all $h_q$. However, for the quadratic case the convergence rate is optimal only for finer quadrature resolution and degrades for coarser quadrature resolution. This result emphasises the importance of subrefining the cut FE cells into sufficiently fine quadrature grid to improve the approximation of the curved boundary. 

%The following experiment aims at evaluating the effect of subrefinement for quadrature.
%The dependence of the solution error measured as the relative $L_2$-norm, i.e. $\|u-u_h\|_{L_2}/\|u_h\|_{L_2}$, and the relative $H^1$-seminorm, i.e. $|u-u_h|_{H^1}/|u_h|_{H^1}$, on the grid size $h_f$ is presented first only for uniform mesh refinements in Figure~\ref{fig:sphere_convergence_h}.\todo[inline]{Pavel: what exactly does $h_q$ variable mean in that plot?}
%We experiment here with subrefinement of all cut elements to improve the quadrature.
%One can see that while the improved quadrature has little effect for linear basis functions,
%it is the key to achieving the optimal convergence rate for quadratic polynomials,
%where the piecewise linear reconstruction of the curved boundary is not sufficient for achieivng the higher order for $h_q = h_f$.

\subsection{Adaptive refinement of the FE grid \label{sec:adaptive_cube}}
\begin{figure*}[]
	\centering
    	\subfloat[][Decomposition of the FE grid into eight subdomains]{
	\includegraphics[width=0.27\textwidth]{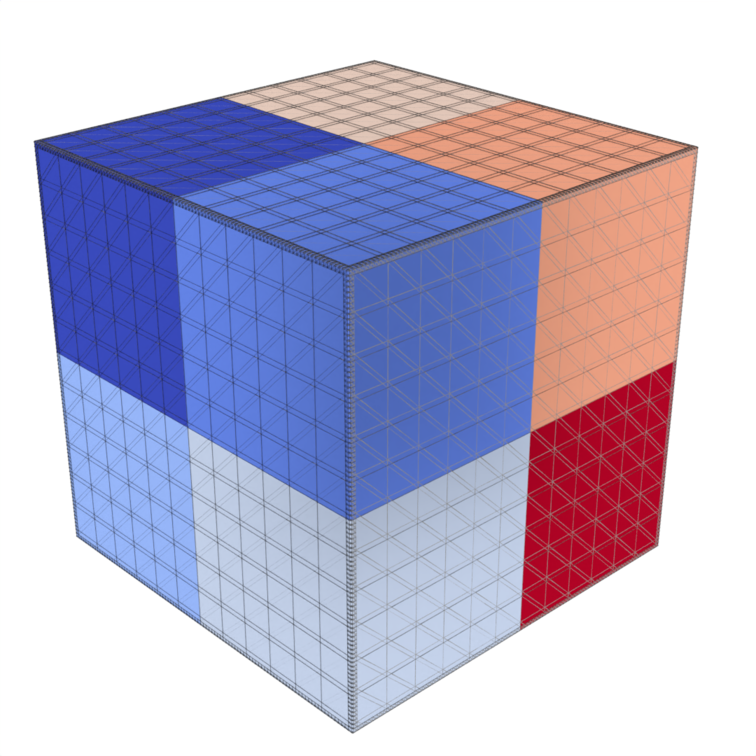} \hspace{0.02\textwidth}
	\includegraphics[width=0.27\textwidth]{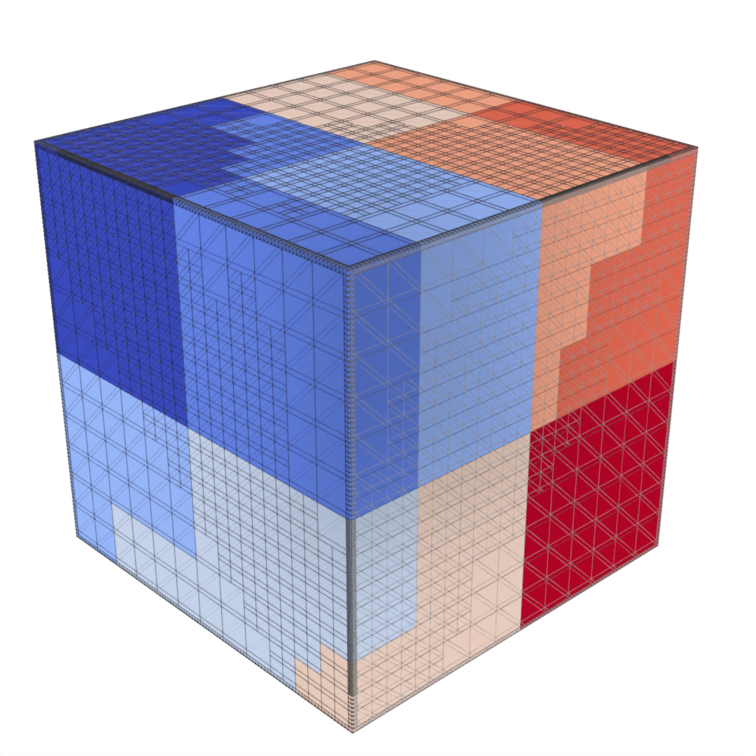} \hspace{0.02\textwidth}
	\includegraphics[width=0.27\textwidth]{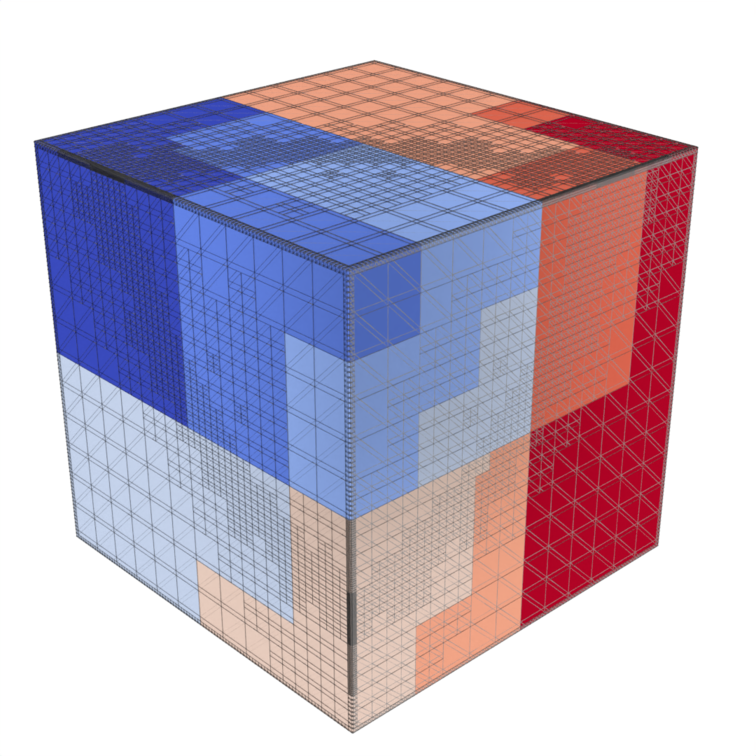}} \\
 	\subfloat[][FE solution]{
	\includegraphics[width=0.27\textwidth]{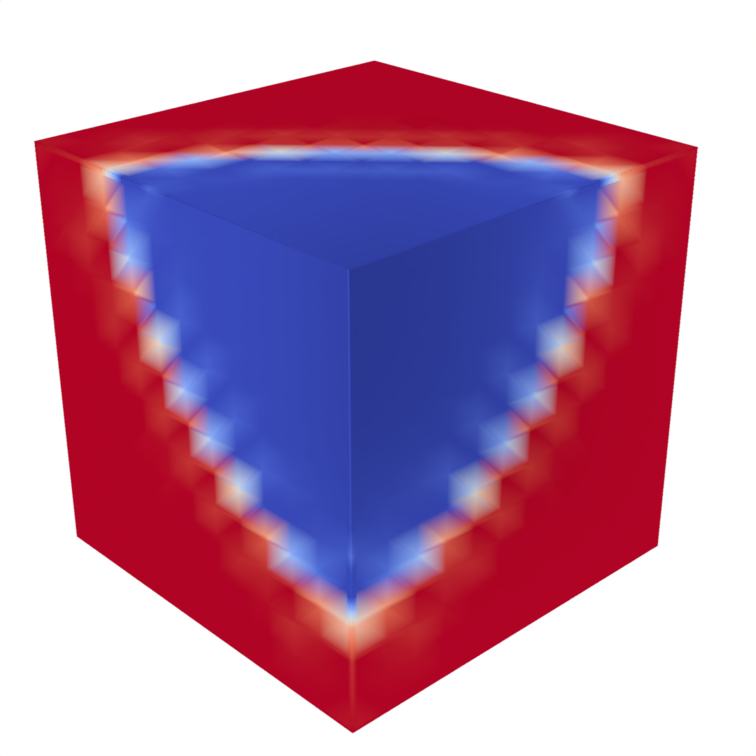} \hspace{0.02\textwidth}
	\includegraphics[width=0.27\textwidth]{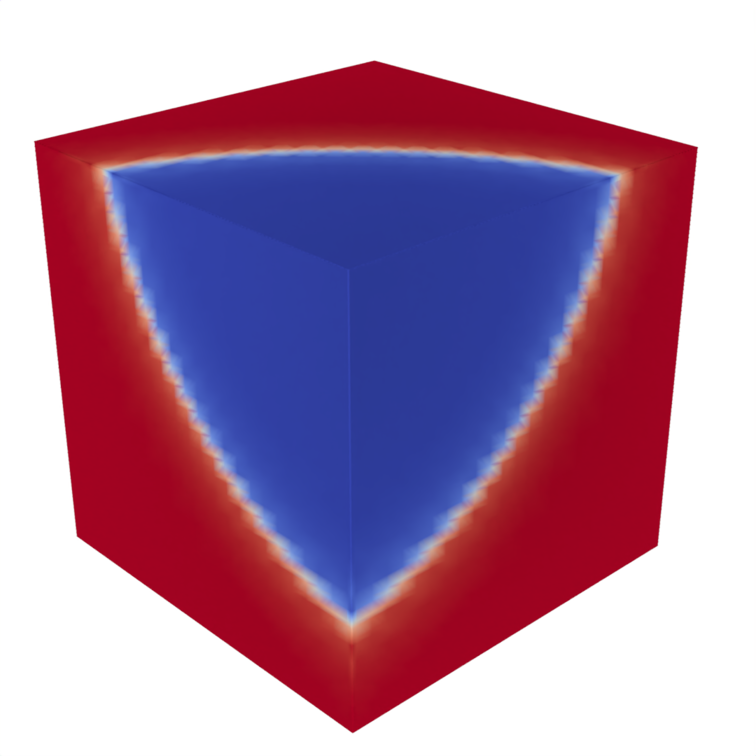} \hspace{0.02\textwidth}
	\includegraphics[width=0.27\textwidth]{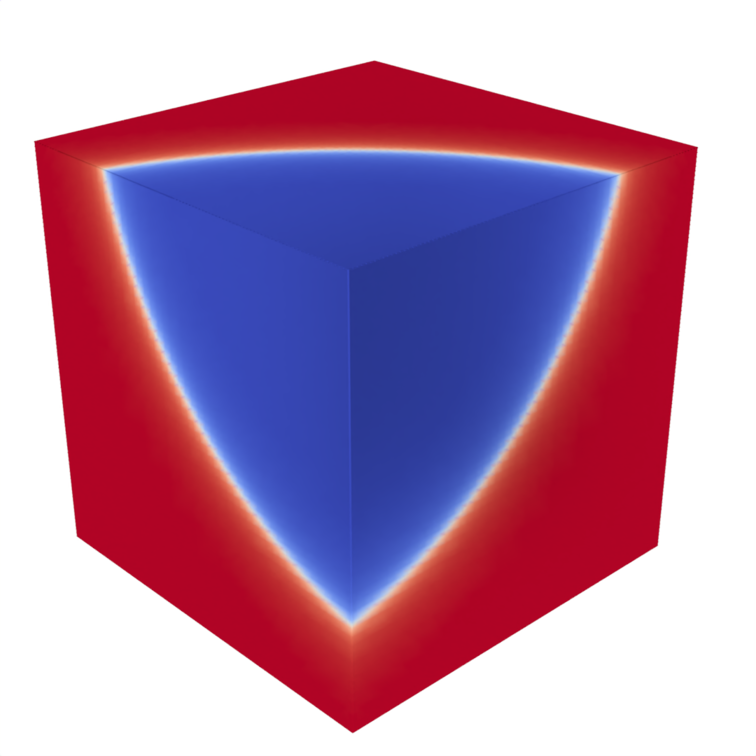} }
	\caption{Adaptively refined grids for the internal layer problem. In each column the FE grid and the respective solution are shown.
	\label{fig:layer_solution}}
\end{figure*}

%--------------------------------------------------------------------------------
%
As the second example we consider an internal layer problem in a unit cube $\Omega = [0,1]^3$ with the prescribed solution \mbox{$u = \arctan(60(r - \pi /3))$}. Here $r = (x_1^2 + x_2^2 + x_3^2)^{1/2}$ is the distance from the origin. Note that the cube domain $\Omega$ is embedded in a larger bounding box $\Omega_\Box$ which generally is not aligned with the unit cube. In this example we test a sequence of uniform and error-driven adaptive FE grid refinements. The coarsest FE grid shown in Figure~\ref{fig:layer_solution} is obtained using 4 uniform refinements of the bounding domain $\Omega_\Box$ which is adaptively refined once and twice, see Figure~\ref{fig:layer_solution}. It is also emphasised that in this example we perform 3 octree refinements of the extraordinary cut FE cells detected using the sharp feature indicator (\ref{eq:sharpCrit}), i.e., $h_q = h_f / 8$ for FE cells containing the corners and edges of $\Omega$. The computation was performed using 1024 cores of the \emph{Salomon} supercomputer at the IT4Innovations National Supercomputing Centre in Ostrava, Czech Republic. Its computational nodes are equipped with two 12-core Intel Xeon E5-2680v3, 2.5 GHz processors, and 128 GB RAM.

\begin{figure}[]
	\centering
	\includegraphics[width=0.48\textwidth]{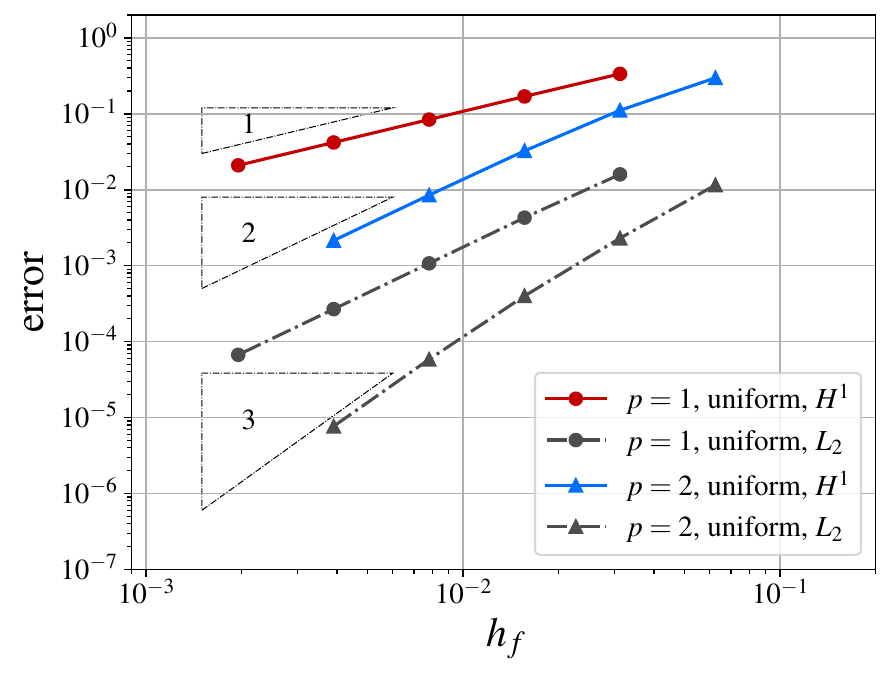} 
	\caption{Dependence of the solution error on $h_f$ for linear ($p=1$) and quadratic ($p=2$) elements for the internal layer problem.}
	\label{fig:convergence_h}
\end{figure}

\begin{figure*}[]
	\centering
	\subfloat[][Relative $L_2$ norm error]{
		\includegraphics[width=0.47\textwidth]{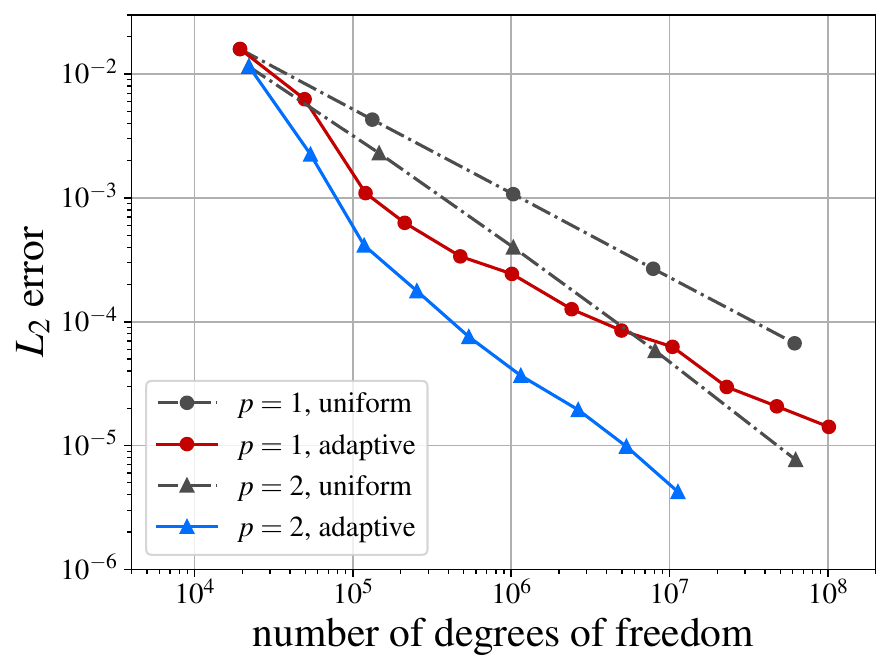} 
	}
	\hspace{0.02\textwidth}%
	\subfloat[][Relative $H^1$-seminorm error]{
		\includegraphics[width=0.47\textwidth]{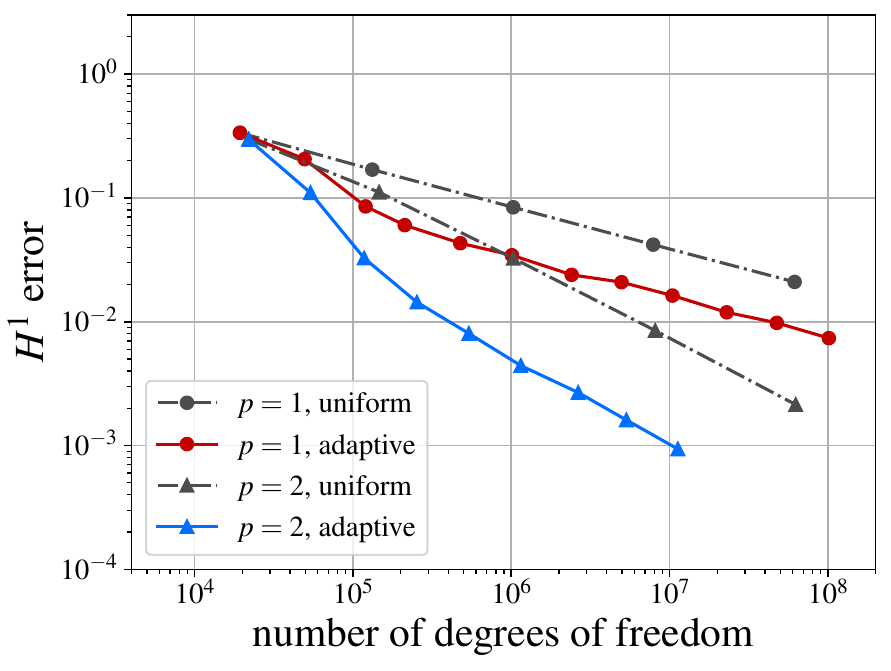}	
	}
	\caption{Dependence of the relative finite element errors on the number of degrees of freedom for uniform and adaptive refinements for the internal layer problem.}
	\label{fig:convergence}
\end{figure*}

The dependence of the solution error measured as the relative $L_2$-norm and the relative $H^1$-seminorm with respect to the FE grid resolution $h_f$ for uniform grid refinements is presented in Figure~\ref{fig:convergence_h}. Here we can see that for this domain with straight faces, we are able to achieve the optimal convergence rates for linear and quadratic basis functions even without using the improved quadrature, i.e., with $h_q = h_f$.
%\todo{Pavel: Is it true for all embedings of our cube? Or does it hold just when the physical and computational cubes are aligned?}
Figure~\ref{fig:convergence} shows the relative $L_2$-norm and the relative $H^1$-seminorm errors with respect to the number of degrees of freedom for both uniform and adaptive refinement. In the adaptive, or error-driven, approach, we refine the elements with the largest $H^1$-seminorm error within each adaptive step such that
approximately 15\% of elements of the whole box are refined. A detailed description of a parallel implementation of this refinement strategy is provided in~\cite{Kus-2017-CND}. It can be observed from Figure~\ref{fig:convergence} that the adaptive refinement achieves optimal convergence rate for both linear and quadratic basis functions with lower relative errors than the uniform refinement. 
In other words it requires more than ten times less degrees of freedom to achieve the same precision as the uniform grid.

%normalised solution errors on the number of degrees of freedom for uniform and error-driven refinements is presented separately for the relative $L_2$-norm and the relative $H^1$-seminorm in Figure~\ref{fig:convergence}.
%In the error-driven approach, we refine the elements with the largest error in the $H^1$-seminorm within each adaptive step such that
%approximately 15 percent of elements of the whole box are refined.
%A detailed description of a parallel implementation of this refinement strategy is provided in~\cite{Kus-2017-CND}.

%--------------------------------------------------------------------------------
\subsection{Strong scalability\label{sec:scalingTest}}
\begin{figure*}[]
	\centering
	\subfloat[][Decomposition of the FE grid]{
		\includegraphics[width=0.3\textwidth]{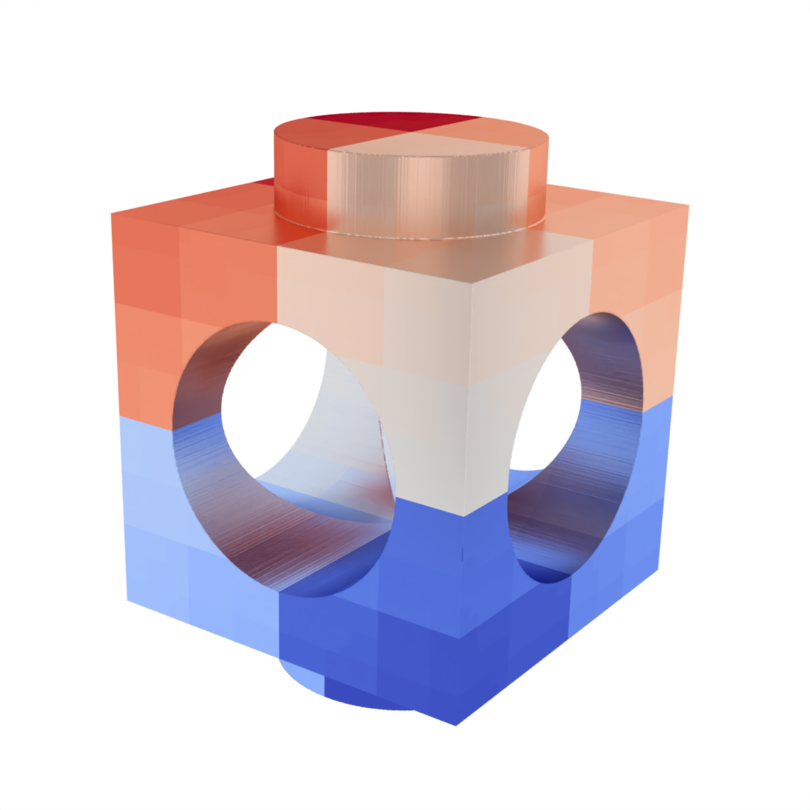} 
	}
	\hspace{0.02\textwidth}%
	\subfloat[][FE solution]{
		\includegraphics[width=0.3\textwidth]{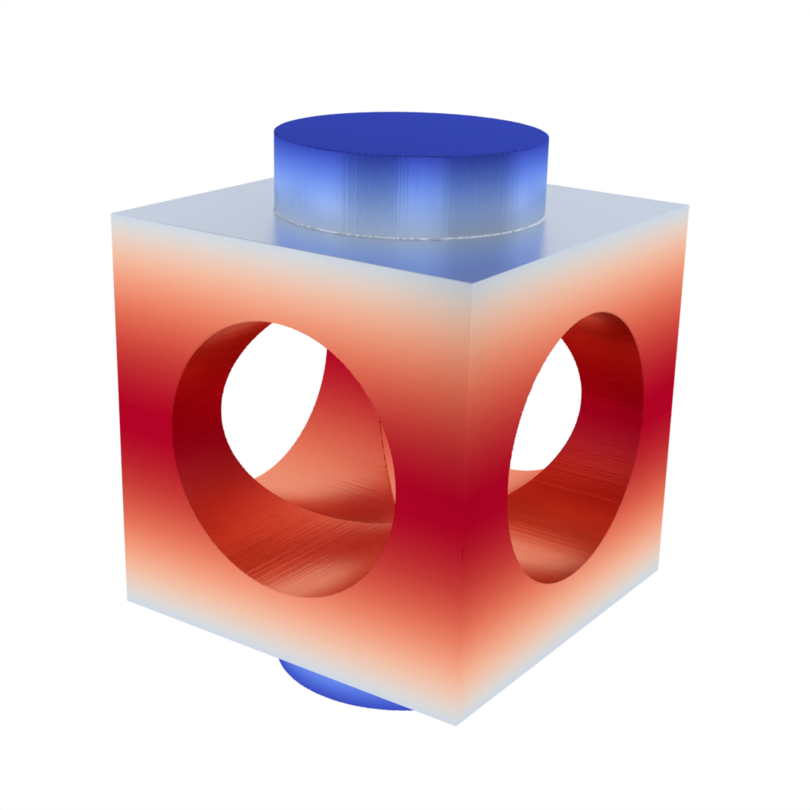}	
	}
	\caption{The handle block geometry. The decomposition of the FE grid into 128 subdomains and the FE solution. }
	\label{fig:handleblock}
\end{figure*}

%--------------------------------------------------------------------------------
%
We assess the strong parallel scalability of our solver by solving the Poisson-Dirichlet problem with a prescribed solution \mbox{$u = \cos(x_3)$}. We consider in this example a handle block geometry as shown in Figure~\ref{fig:handleblock}. In strong scaling the finite element size is fixed and the number of processors is continuously increased. The signed distance function of the handle block is obtained using constructive solid geometry (CSG)~\cite{Ricci1973} of three cylinders and a box. In this example we use a base grid obtained by refining the bounding box $\Omega_\Box$ uniformly six times, where grid cells far from the boundary are coarsened by default. We obtain from the base grid an FE grid through five levels of refinements towards the boundary with one additional refinement level for quadrature of extraordinary cells. The resulting FE grid contains 12 million trilinear elements and 11 million degrees of freedom. 
%\todo{not clear here the sizing of geometry grid and quadrature grid}
%We solve the Poisson problem with a prescribed solution (\ref{eq:cosine}).
%The right-hand side vector and boundary conditions are set according to the solution. 
%The FE grid is obtained by 6 uniform refinements, 3 unrefinements, 5 boundary refinements and 1 octree subrefinement of extraordinary cells.
%The resulting problem contains 12 million elements and 11 million DOFs.

The strong scaling test was performed on the \emph{Salomon} supercomputer at the IT4Innovations National Supercomputing Centre in Ostrava, Czech Republic.
Its computational nodes are equipped with two 12-core Intel Xeon E5-2680v3, 2.5 GHz processors, and 128 GB RAM.
The number of subdomains, i.e., the number of processes in the parallel computation, ranges from 128 to 1024.

\begin{figure*}[]
	\centering
	\subfloat[][\label{fig:scaling_handleblock}]{
		\includegraphics[width=0.47\textwidth]{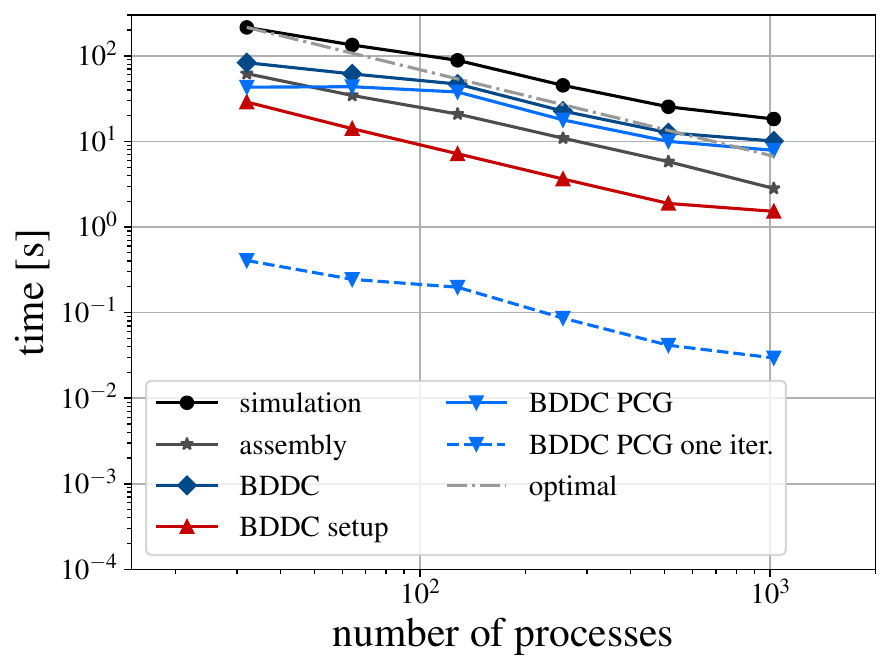} 
	}
	\hspace{0.02\textwidth}%
	\subfloat[][\label{fig:pcg_handleblock}]{
		\includegraphics[width=0.47\textwidth]{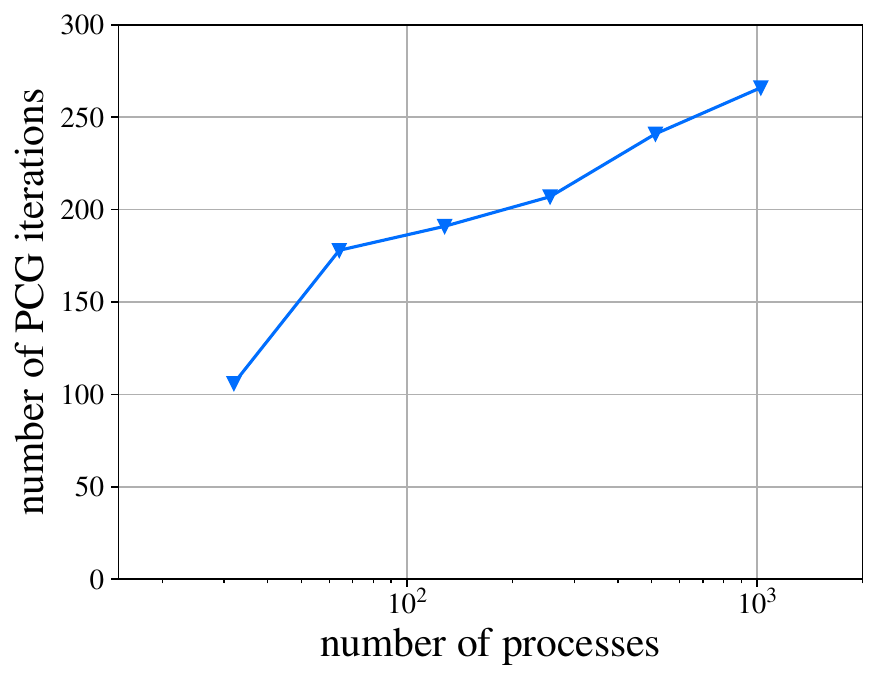}	
	}
    \caption{Strong scaling test on the handle block geometry, 6 uniform refinements, 3 unrefinements, 5 refinements towards boundary, 1 octree refinement of extraordinary cells, finite element order 1, 3-level BDDC. Problem with 12 million elements and 11 million DOFs. Timing of the important parts of the solver (a) and number of PCG iterations (b).}
	\label{fig:strong_scaling}
\end{figure*}

The results of the strong scaling test are presented in Figure~\ref{fig:scaling_handleblock}, where the runtimes of the important components of the solver are analysed separately. These include assembly of the matrices and time for the BDDC solution. The latter is further divided into the time spent in the BDDC setup and in PCG iterations. In addition, we report the total time necessary for the whole simulation.
Optimal strong scalability corresponds to halving the computational time every time the number of subdomains is doubled. In the logarithmic scale, this corresponds to a straight line, marked as `optimal' in Figure~\ref{fig:scaling_handleblock}.
It is evident from the figure that the scalability of assembly and BDDC setup is optimal. However, the PCG iterations within BDDC do not scale optimally, which leads to a suboptimal scaling of the overall simulation. The suboptimal scaling in the PCG can be attributed to the growing number of iterations as identified in Figure~\ref{fig:pcg_handleblock}. 
Consequently, we also plot the time for one PCG iteration showing the favourable strong scalability of our implementation.
%the whole simulation is slightly suboptimal. If we trace the cause of this, it is the BDDC solver, namely the part spent in the PCG iterations that is not scaling in an optimal speed.
%However, this increase can be atributed to the growing number of PCG iterations (see the right plot in Figure~\ref{fig:pcg_handleblock}),
%and we consider it to be acceptable.

%--------------------------------------------------------------------------------
\subsection{Application to complex engineering CAD models \label{sec:engCAD}}
%--------------------------------------------------------------------------------
%
%\todo[inline]{Matt's comment: We should say that the STL files are obtained by meshing CAD, and maybe cite opencascade or whatever you used, along with the meshing parameters used. (BTW, have you guys experimented with chatgpt? Ask it to generate OpenCascade code to visit the faces of a brep, or sample normals on a brep face, or generate cmake files to build it. ChatGPT is great for figuring out annoying apis)}

\begin{figure*}[]
 	\centering
 	\includegraphics[width=0.8\textwidth]{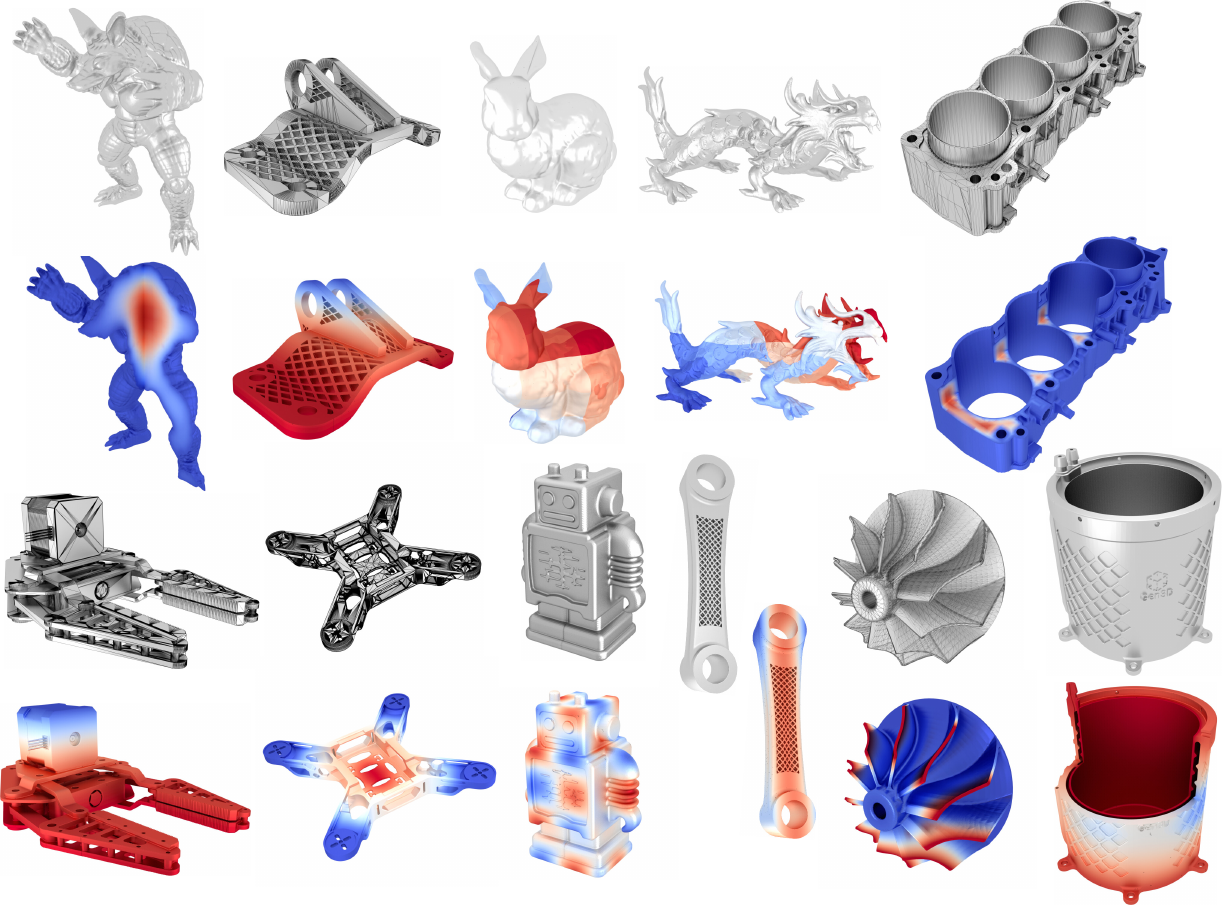} 
 	\caption{Complex geometries used for the robustness test (from top left to bottom right): \emph{armadillo}, \emph{bracket}, \emph{bunny}, \emph{dragon},  \emph{engine}, \emph{gripper}, \emph{quadcopter}, \emph{robo}, \emph{spanner}, \emph{turbo}, and \emph{cover}.}
 	\label{fig:all-geometry}
 \end{figure*}

In this section we first assess the robustness of the proposed method for solving Poisson-Dirichlet problems on various complex geometries from computer graphics and engineering, see Figure~\ref{fig:all-geometry}. The engineering models are given in STEP format and are converted into sufficiently fine STL meshes using FreeCAD. The computer graphics models are given as STL meshes. The STL meshes are immersed in the geometry grid with resolution smaller than the resolution of the quadrature grid, i.e., $h_g \leq h_q$. For all geometries, the base grid is obtained by five uniform refinements of $\Omega_\Box$. The FE grid is obtained by six refinements of the base grid towards the boundary. The quadrature grid is obtained by refining the extraordinary cut FE cells once. 
We consider the prescribed solution $u = \cos(x_3)$, trilinear Lagrange polynomials on FE cells as the basis, and three levels in the BDDC method.

The computations are performed on the \emph{Karolina} supercomputer at the IT4Innovations National Supercomputing Centre in Ostrava, Czech Republic. The computational nodes are equipped with two 64-core AMD 7H12 2.6 GHz processors, and 256 GB RAM.

We decompose the FE grid into 256 subdomains assigned to the same number of computer cores. Table~\ref{tab:problems_overview} shows the number of degrees of freedom, number of elements, number of iterations to reach the relative residual precision of $10^{-6}$, and computational runtimes for each geometry. In addition, we compare the runtimes of the domain decomposition solver with the parallel sparse direct solver MUMPS~\cite{Amestoy:2019}.
It is evident from Table~\ref{tab:problems_overview} that the BDDC solver is reasonably robust and able to converge for all these geometries requiring from 257 to 913 PCG iterations.
While MUMPS is also able to provide the solutions for all geometries when provided enough memory, it is consistently slower than the BDDC solver by a factor ranging from 1.2 to 6.3 depending on the problem.

\begin{table*}
	\centering
\begin{tabular}{l | r r | r r | r | r r r | r}
%# problem nrefu nuref nrefb ndofsextr ndofi nelementsextr iter tlevelset trefinements tassembly tsetup tkrylov tloading tbddcml twhole twholemumps
            &  \multicolumn{2}{c|}{DOFs} & &      & \multicolumn{5}{c}{Time [s]} \\
Problem     &  Global & Interf. & Elems. & Iters. & Assembly & Setup & PCG   & Solve BDDC & Solve MUMPS \\
\hline
Armadillo   &  9.2M   & 384k    & 11.9M  & 360    &  20.1 &  2.9  &  72.1 & 80.0   & 238.5  \\ 
Bracket     &  11.3M  & 523k    & 14.4M  & 343    &  28.9 &  3.3  &  77.0 & 87.9   & 275.8  \\
Bunny       &  12.8M  & 463k    & 16.7M  & 450    &  29.5 &  4.3  & 110.5 & 122.4  & 359.8  \\
Dragon      &  5.2M   & 280k    & 6.8M   & 375    &  11.2 &  1.7  &  79.3 & 83.6   & 125.3  \\
Engine      &  11.6M  & 533k    & 14.9M  & 436    &  32.7 &  4.0  &  89.3 & 105.4  & 290.9  \\
Gripper     &  10.1M  & 443k    & 12.7M  & 913    &  25.9 &  3.2  & 160.1 & 173.9  & 230.9  \\
Quadcopter  &  4.5M   & 277k    & 5.6M   & 264    &  12.0 &  1.5  &  59.5 & 64.0   &  95.9  \\
Robo        &  11.2M  & 378k    & 14.6M  & 257    &  27.5 &  3.9  &  65.2 & 77.6   & 260.2  \\
Spanner     &  5.2M   & 338k    & 6.7M   & 503    &  12.4 &  1.9  &  97.2 & 102.1  & 124.0  \\
Turbo       &  2.0M   & 620k    & 2.5M   & 270    &  50.6 &  6.1  &  58.5 & 79.0   & 498.2  \\
Cover&  11.6M  & 424k    & 15.0M  & 209    &  35.1 &  10.1 &  33.6 & 45.0   & 459.7 
\end{tabular}                 
\caption{\label{tab:problems_overview}Robustness test for different geometries, 5 uniform refinements, 3 unrefinements, 6 refinements towards boundary, 1 refinement for quadrature.
In the legend,
`Interf.' is \emph{interface}, `Iters.' is \emph{number of iterations}, `Setup' is \emph{BDDC setup}, `Solve BDDC' is the overall \emph{BDDC solver time} and `Solve MUMPS' is the overall time using the MUMPS distributed sparse direct solver.}
\end{table*}

\subsection{Weak scalability studies \label{sec:weakScalability}}
\begin{figure*}[]
	\centering
	\subfloat[][\label{fig:scaling_emotor}]{
		\includegraphics[width=0.47\linewidth]{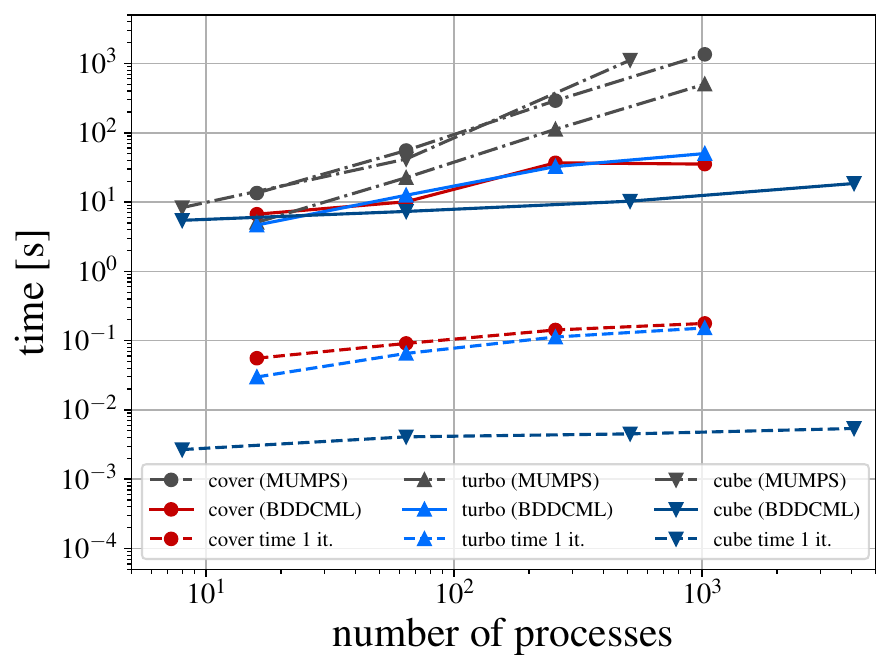} 
	}
	\hspace{0.02\textwidth}%
	\subfloat[][\label{fig:pcg_emotor}]{
		\includegraphics[width=0.47\linewidth]{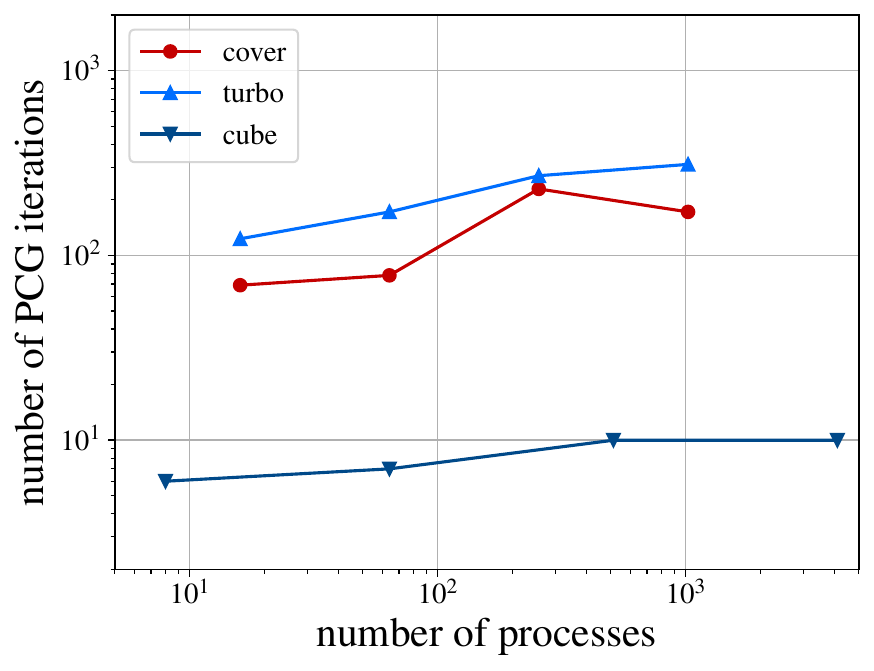}	
	}
	\caption{Weak scaling of the BDDC and the direct solver (MUMPS) for the \emph{cover} and the \emph{turbo} problems (left), and the corresponding numbers of iterations required by the BDDC method to reach relative residual precision of $10^{-6}$ (right).}
	\label{fig:emotor_weak_scaling}
\end{figure*}

%--------------------------------------------------------------------------------
%
Finally, we test the weak scalability of the proposed approach first using  the \emph{turbo} and \emph{cover} geometries shown in Figure~\ref{fig:all-geometry}.
In the left part of Figure~\ref{fig:emotor_weak_scaling},
we present computational time for an increasing number of processes,
and we compare the three-level BDDC method with using a distributed sparse Cholesky factorization by MUMPS. During the test, the number of processes increases from 16 to 1\,024, and we perform one refinement of elements toward boundary within each step.
The problem size grows approximately four times with each of these refinements, starting at 660 thousand unknowns and finishing with 46 million unknowns in the global system for the \emph{cover} problem, 
and from 290 thousand to 20 million unknowns for the \emph{turbo} problem.
For both problems, the local problem size is kept approximately fixed, around 45 thousand unknows per subdomain for the cover problem, and around 19 thousand for the turbo problem.

We can clearly observe the growing time of the MUMPS solver corresponding to the increasing complexity of the sparse Cholesky factorization of the distributed matrix. 
On the other hand, the three-level BDDC method requires significantly less time to reach the desired precision,
although the solution time grows as well.
Most of this growth can be attributed to the increasing number of iterations, which is plotted separately in the right part of Figure~\ref{fig:emotor_weak_scaling}.
For the largest presented problem, the three-level BDDC method is 38 times faster than the MUMPS solver.

Furthermore, we consider the unit cube geometry introduced in Section~\ref{sec:adaptive_cube} to evaluate the weak scalability of our approach for a more structured problem.
We can see from Figure~\ref{fig:emotor_weak_scaling} that the weak scalability for this problem is almost optimal with the number of iterations not growing beyond 10.

%\begin{figure}[htp]
%	\centering
%		\includegraphics[width=0.8\linewidth]{./figs/plots/scaling_emotor} 
%		\includegraphics[width=0.8\linewidth]{./figs/plots/scaling_emotor_iterations}	
%	\caption{Weak scaling of the BDDC and the direct solver (MUMPS) for the \emph{cover} and the \emph{turbo} problems (left), and the corresponding numbers of iterations required by the BDDC method to reach relative residual precision of $10^{-6}$ (right).}
%	\label{fig:emotor_weak_scaling}
%\end{figure}

%\todo[inline]{Should we add also the engine block?}
%\todo[inline]{Pavel: Given the focus of this paper, I think that having more real-world challenging geometries would be a plus}

	%--------------------------------------------------------------------------------
\section{Conclusions}
\label{sec:conclusions}
%--------------------------------------------------------------------------------

%\begin{itemize}
%	\item[-] neural nets for representing and processing of signed distance functions~\cite{}
%\end{itemize}

% Conclude the geometry grid resolves the domain and bounds the geometric error
We introduced an immersed finite method that uses three non-boundary conforming background grids to solve the Poisson problem over complex CAD models. The three grids, namely the geometry, finite element, and quadrature grids, are each constructed to represent the geometry, discretise the physical field, and for integration, respectively. In our examples, we use an implicit description of the geometries in the form of a scalar-valued signed distance function. This is obtained by evaluating the distance from the grid points to the sufficiently fine triangle surface mesh (STL mesh) approximating the exact CAD surface. To capture geometrical features up to a user-prescribed resolution, the geometry grid is refined within a narrow band from the boundary. In our computations, we choose the geometry grid to have the highest resolution among the three grids. Our examples have shown that the resolution of the geometry grid plays a crucial role in determining the boundary approximation error as well as limiting the integration error.

% Conclude all the treatments yield optimal FE convergence
% Higher order integration is performed through the generation of quadrature grid. Mention the limitation: in practice, the refinement for quadrature grid is limited to two due to the escalating integration cost. 
The finite element grid for discretisation is constructed by considering the features of the physical field. In cases where the solution features coincide with the boundary, we have presented a technique to correctly identify FE grid cells that are cut by the domain boundary, i.e., cut cells, including those containing sharp and other small geometrical features. In the integration step, such cut cells are integrated by first constructing the quadrature grid using a bottom-up strategy. The leaf cells of the quadrature grid are subsequently tessellated using marching tetrahedra. We have demonstrated with an internal layer example that an optimal convergence rate can be obtained for both linear and high-order polynomial basis functions using two levels of refinement for the quadrature grid. When FE basis functions cover only a tiny part of the physical domain, we use a cut-cell stabilisation technique that extends the basis support over the nearest active cells. 

% Conclude the parallelisation enables computing complex geometries
% 1. Importance of weights in the partitioning for p4est
% 2. Solver has to handle multi-component subdomains
% 3. The candidate for basis extension (stabilisation) is currently limited to the subdomains
%The FE grid is subsequently partitioned and distributed across several processors. 

% Future works
There are several promising extensions of the proposed three-grid immersed finite element method. One potential future direction of this research is to use a neural network for representing the implicit signed distance function, as demonstrated in~\cite{park2019deepsdf} for point cloud input data. The use of point clouds as geometric input is gaining interest due to its direct connection with inspection and monitoring of engineering products, as shown in~\cite{Balu2023}. 
The application of the current approach beyond the Poisson equation to other second-order elliptic equations, like linear elasticity, is straightforward, as the only significant change involves the integrands in the weak form. However, the consideration of other kinds of linear and nonlinear equations, like elastodynamics, requires further analysis. Additionally, the proposed cut-cell identification algorithm can be modified to detect long-sliver cracks and other features if they are considered necessary in some applications.  One possible strategy is to use bottom-up adaptive mesh refinement approach for the FE grid, starting from the finest resolution of the geometry grid. Lastly, another promising extension is the refinement of the finite element grid according to an a posteriori error indicator instead of relying on a priori knowledge of the solution features, as done in the current work.

    \section*{Acknowledgments}
     JS and PK have been supported by the Czech Science Foundation through GAČR 23-06159S, and by the Czech Academy of Sciences through RVO:67985840. Computational time on the Salomon and Karolina supercomputers has been provided thanks to the support of the Ministry of Education, Youth and Sports of the Czech Republic through the e-INFRA CZ (ID:90140). EF has been supported by the Royal Society Research Grants through~RGS{\textbackslash}R2{\textbackslash}222318. 
        
     %\appendix
	
	%\input{appendix}

	%-------------------------------------------------------------------------------
	\bibliographystyle{elsarticle-num-names}
	\bibliography{immersedParallel.bib}
	
\end{document}